\documentclass[11pt]{amsart}

\usepackage{comment}
\usepackage{amsmath}							
\usepackage{amssymb}
\usepackage{amsthm}
\usepackage{amscd}
\usepackage{amsfonts}
\usepackage{mathabx}
\usepackage[all]{xy}
\usepackage{colonequals}
\usepackage{upgreek}
\usepackage{graphicx}
\usepackage{etoolbox}

\makeatletter
\DeclareRobustCommand{\cev}[1]{%
  \mathpalette\do@cev{#1}%
}
\newcommand{\do@cev}[2]{%
  \fix@cev{#1}{+}%
  \reflectbox{$\m@th#1\vec{\reflectbox{$\fix@cev{#1}{-}\m@th#1#2\fix@cev{#1}{+}$}}$}%
  \fix@cev{#1}{-}%
}
\newcommand{\fix@cev}[2]{%
  \ifx#1\displaystyle
    \mkern#23mu
  \else
    \ifx#1\textstyle
      \mkern#23mu
    \else
      \ifx#1\scriptstyle
        \mkern#22mu
      \else
        \mkern#22mu
      \fi
    \fi
  \fi
}
\makeatother

\usepackage{caption}
\usepackage{subcaption}

\usepackage{euler}

\usepackage{extarrows}

\usepackage[colorlinks, linktocpage, citecolor = purple, linkcolor = blue]{hyperref}
\usepackage[noabbrev,capitalize]{cleveref}
\usepackage{color}
\usepackage{faktor}

\usepackage{tikz}	
\usetikzlibrary{calc}
\usepackage{float} 
\usetikzlibrary{matrix}
\usetikzlibrary{patterns}
\usetikzlibrary{positioning}
\usetikzlibrary{decorations.pathmorphing}
\usetikzlibrary{cd}
\definecolor{lilla}{RGB}{199,85,255}
\usepackage{ytableau}
\usepackage{fullpage}
\usepackage[shortlabels]{enumitem}

\usepackage{microtype}

\newtheorem{maintheorem}{Theorem}	
\newtheorem{theorem}{Theorem}[section]
\newtheorem{lemma}[theorem]{Lemma}
\newtheorem*{lemma*}{Lemma}
\newtheorem{proposition}[theorem]{Proposition}
\newtheorem{corollary}[theorem]{Corollary} 
\newtheorem*{corollary*}{Corollary} 

\newtheorem*{conjecture*}{Conjecture}

\theoremstyle{definition}
								
\newtheorem{definition}[theorem]{Definition}
\newtheorem{example}[theorem]{Example}
\newtheorem*{example*}{Example}

\AtEndEnvironment{construction}{\qed}

\newtheorem{claim}[theorem]{Claim}
\crefname{claim}{Claim}{Claims}

\newtheorem*{remark*}{Remark}
\newtheorem{question}[theorem]{Question}
\newtheorem{question*}{Question}

\newcommand{\casetitle}[2][]{%
  \par\medskip
  \noindent\textbf{Case (#2)}%
  \if\relax\detokenize{#1}\relax\else: \textit{#1}\fi.\quad
}

\makeatletter
\def\@tocline#1#2#3#4#5#6#7{\relax
	\ifnum #1>\c@tocdepth % then omit
	\else
	\par \addpenalty\@secpenalty\addvspace{#2}%
	\begingroup \hyphenpenalty\@M
	\@ifempty{#4}{%
		\@tempdima\csname r@tocindent\number#1\endcsname\relax
	}{%
		\@tempdima#4\relax
	}%
	\parindent\z@ \leftskip#3\relax \advance\leftskip\@tempdima\relax
	\rightskip\@pnumwidth plus4em \parfillskip-\@pnumwidth
	#5\leavevmode\hskip-\@tempdima
	\ifcase #1
	\or\or \hskip 1em \or \hskip 2em \else \hskip 3em \fi%
	#6\nobreak\relax
	\dotfill\hbox to\@pnumwidth{\@tocpagenum{#7}}\par
	\nobreak
	\endgroup
	\fi}
\makeatother

\newcommand{\PP}{\mathbb{P}}

\DeclareMathOperator{\gon}{gon}

\let\ddiv\relax
\DeclareMathOperator{\ddiv}{div}

\makeatletter
\newcommand{\euleroperator}[1]{%
\mathop{\text{\usefont{OT1}{cmr}{m}{n}#1}}\nolimits}
\makeatother

\let\deg\relax
\newcommand{\deg}{\euleroperator{deg}}
\let\gon\relax
\newcommand{\gon}{\euleroperator{gon}}

\title{On gonality-tight graphs}

\author{\v Simun Dropulji\' c}
\address{Mathematical Institute, University of St Andrews, St Andrews KY16 9SS, UK}
\email{\href{mailto:sd304@st-andrews.ac.uk}{sd304@st-andrews.ac.uk}}

  \author{Yoav Len}
   \address{Mathematical Institute, University of St Andrews, St Andrews KY16 9SS, UK}
 \email{\href{mailto:yoav.len@st-andrews.ac.uk}{yoav.len@st-andrews.ac.uk}}

\begin{document}

\begin{abstract}

 We address a question posed by Fessler--Jensen--Kelsey--Owen regarding graphs whose second gonality is greater than the first by exactly $1$.  We answer the question affirmatively under a 
  stronger condition, thereby characterising the entire gonality sequence for a large family of graphs, that we refer to as quasi-banana graphs. We prove a structure theorem for those graphs, and show that they are all obtained via an inductive process by gluing together complete and banana graphs under certain rules.

%We provide computational evidence in support of a positive answer in general and suggest a stronger version of the original question. 
\end{abstract}

\maketitle

\setcounter{tocdepth}{2}
\tableofcontents

%%%%%%%%%%%%%%%%%%%%%%%%%%%%%%%%%%%%%%%%%%%%%%%%%%%%%%

\section{Introduction}\label{sec:introduction} 

The $k$-th gonality of a graph, denoted $\gon_k(G)$, is the smallest degree of a divisor of rank $k$. The gonality sequence $\gon_1(G), \gon_2(G),\ldots$, controls much of its geometric properties.  For instance, 
if any element $\gon_k(G)$ equals $2k$ for $k$  strictly smaller than the genus, then the graph is hyperelliptic \cite{Coppens_Clifford, Len_Clifford}. Questions about the gonality of graphs first appeared in the seminal paper of Baker and Norine in which they introduced the rank of divisors and discussed the notion of hyperelliptic graphs \cite{BakerNorine2007}. The connections between the gonality and covers of trees were then studied in  \cite{CH, DraismaVargas_Catalan, MeloZheng_trigonal} and their moduli space in \cite{CoolsDraisma_gonalityModuli}.
The first reference to the gonality sequence of a graph as a whole appeared in \cite{CoolsPanizzut_gonalityComplete}, where its values for the complete graphs were determined.  

The gonalities of algebraic curves have been studied since the early days of algebraic geometry, but, as far as we can tell, the term ``gonality sequence" first appeared in \cite{CoppensKeemMartens_primitiveLinearSerie} and subsequently in \cite{LangeMartens_gonalitySequence, LangeNewstead_CliffordIndices}. For the general curve, its gonality sequence is completely determined by the Brill--Noether theorem, and the gonality sequence of hyperelliptic curves is known  as well  \cite[page 114]{KatoMartens_gonalityWithInvolution}.
However, there are still many open questions regarding the possible gonality sequences of special curves, and the interplay between tropical and algebraic geometry has the potential for new results. For instance, the authors of \cite{CDJP_BrillNoetherForP1xP1}
 use the gonality sequence of bipartite graphs to compute those of curves in 
 $\PP^1\times\PP^1$.

It is natural to ask which sequences may appear as the gonality sequence of a curve or a graph, and to describe the locus of objects with this sequence (see  \cite[Question 1.3]{ADMYY_GonalitySequencesOfGraphs}).  As shown in \cite[Theorem 1.4]{ADMYY_GonalitySequencesOfGraphs}, the gonality sequence is completely determined by the first gonality for graphs of genus smaller than 6. 
The  sequence 
has been partially determined in various cases of generalised banana graphs \cite{CastryckCools_NewtonAndGonalities,FJKO_semigroupOfGonality, ADMYY_GonalitySequencesOfGraphs} and chess graphs
 \cite{Speeter_RookGraphs, CDDKLMS_ChessGraphs, JMSWY_FerresGraphs}, but not much  is known in general.

\begin{figure}
\centering
\begin{tikzpicture}[scale=2.2, node/.style={circle, draw, fill=black, inner sep=1.2pt}]

    \begin{scope}[scale=0.8]
        \node[node] (v) at (0,0) {};
        \node[node] (u11) at (2,0) {};
        \node[node] (u12) at (1.81262,0.84524) {};
        \node[node] (u13) at (1.28558,1.53209) {};
        \node[node] (u14) at (0.51764,1.93185) {};
    \end{scope}
    
    % Edges for first graph
    \draw (u11)--(u12)--(u13)--(u14)--(u11)--(u13)--(u11);
    \draw (u12)--(u14);
    \draw[bend left=6] (v) to (u11);
    \draw[bend right=6] (v) to (u11);
    \draw[bend left=6] (v) to (u12);
    \draw[bend right=6] (v) to (u12);
    \draw[bend left=6] (v) to (u13);
    \draw[bend right=6] (v) to (u13);
    \draw[bend left=6] (v) to (u14);
    \draw[bend right=6] (v) to (u14);

    % Vertices for the second graph
    \node[node] (v1) at (2.2,0.3) {};
    \node[node] (v2) at (3,0.3) {};
    \node[node] (v3) at (3.8,0.3) {};
    \node[node] (v4) at (4.6,0.3) {};
    \node[node] (v5) at (5.4,0.3) {};

    % Edges for the secnd graph
    \draw[bend left=20] (v1) to (v2);
    \draw[bend right=20] (v1) to (v2);
    
    \draw[bend left=25] (v3) to (v2);
    \draw (v3) to (v2);
    \draw[bend right=25] (v3) to (v2);
    
    \draw[bend left=30] (v3) to (v4);
    \draw[bend right=10] (v3) to (v4);
    \draw[bend left=10] (v3) to (v4);
    \draw[bend right=30] (v3) to (v4);

    \draw[bend left=40] (v5) to (v4);
    \draw[bend right=20] (v5) to (v4);
    \draw[bend left=20] (v5) to (v4);
    \draw[bend right=40] (v5) to (v4);
    \draw (v4) to (v5);
\end{tikzpicture}
\caption{Two quasi-banana graphs of gonality $5$.}
\label{fig:gon5example}
\end{figure}
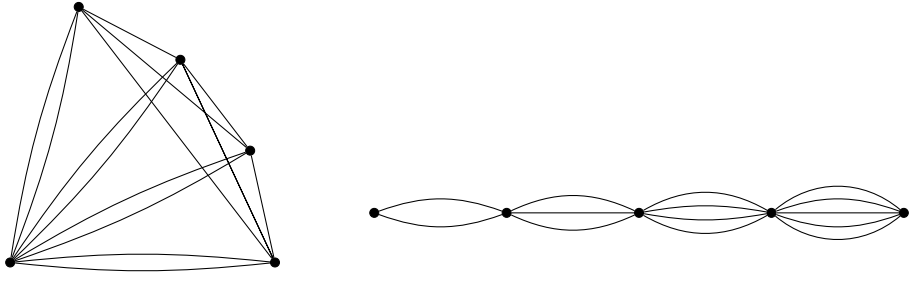

 In the case of algebraic curves, it has recently been shown that 
 much of their geometry is determined whenever the first two elements of the sequence differ by $1$ \cite[Lemma 7.3]{FJKO_semigroupOfGonality}. 
The authors then go on to ask whether  
the same holds for graphs \cite[Question 1.3]{FJKO_semigroupOfGonality}. In this paper, we refer to  graphs or curves whose second gonality exceeds the first by $1$  
as \emph{gonality-tight} (see \cref{def:gonalityTight}).
A key example 
of gonality-tight graphs are the generalised banana graphs 
$B_{a,a}^*$, for which  the second gonality is realised by a divisor with a single chip on each vertex, apart from one, which has two chips \cite[Lemma 5.6]{ADMYY_GonalitySequencesOfGraphs}. 
In the current paper, we examine when the second gonality of a gonality-tight graph is realised by a divisor of this form. As we show, the existence of such a divisor has an effect on the structure of the graph.  In fact, this happens 
exactly for special graphs that we refer to as \emph{quasi-banana graphs} (see \cref{fig:gon5example} for an illustration and \cref{def:quasiBananaGraph} for the exact definition).

\begin{maintheorem}\label{mainTheorem:alternativeDefinitonForQBGraphs}
    A graph $G$ is a quasi-banana graph if and only if it is gonality-tight and there exists a vertex $v\in V(G)$ such that the divisor
    \[
    D=v+\sum_{w\in V(G)}w
    \]
    realises the second gonality.
\end{maintheorem}

We then proceed to examine additional properties of quasi-banana graphs and, in particular,   fully answer \cite[Question 1.3]{FJKO_semigroupOfGonality} for this family of graphs. 

\begin{maintheorem}\label{mainTheorem:propertiesOfQBGraphs}
    Let $G$ be a quasi-banana graph with first gonality $\gon_1(G)=k$. Then the following hold.
    \begin{enumerate}
        \item The genus of $G$ is $g=\binom{k}{2}$.
        \item The gonality sequence of $G$ is given by
        \[
        \gon_r(G)=
        \begin{cases}
            l(k+1)-h & \text{if } r < g, \\
            g+r  & \text{if } r\geq g,
        \end{cases}
        \]  
        where $1\leq l\leq k-2$ and $0\leq h\leq l$ are the uniquely determined integers such that $r=\frac{l(l+3)}{2}-h$.
        \item If $k\neq 2$, all divisors realising the second gonality are linearly equivalent.
        \item Let $D$ be a divisor realising the second gonality. Then, for  every integer $l$ such that $1\leq l \leq k$, the divisor $l\cdot D$ realises the $\frac{l(l+3)}{2}$-th gonality. In particular, $(k-2)\cdot D$ is equivalent to the canonical divisor $K_G$.
    \end{enumerate}    
\end{maintheorem}

Note that quasi-banana graphs that share the same first gonality can still have significantly different structures (as seen in  \cref{fig:gon5example} for instance).

Based on over 600 examples of non-isomorphic gonality-tight graphs obtained computationally (three of which can be seen in \cref{fig:examples}, additional examples are available upon request), we believe that the results hold in greater generality.

\begin{conjecture*}
    Let $G$ be a gonality-tight graph with first gonality $k$ for some integer $k$.
    \begin{enumerate}
        \item The genus of $G$ is equal to $g=\binom{k}{2}$.
        
        \item The gonality sequence of $G$ is given by
        \[
        \gon_r(G)=
        \begin{cases}
            l(k+1)-h & \text{if } r < g,\\
            g+r  & \text{if } r\geq g,
        \end{cases}
        \]  
        where $1\leq l\leq k-2$ and $0\leq h\leq l$ are uniquely determined integers such that $r=\frac{l(l+3)}{2}-h$.

        \item If $D$ realises the second gonality, then $(k-2)\cdot D\sim K_G$. Moreover, for every integer $1\leq l\leq k-2$, the divisors realising the $\frac{l(l+3)}{2}$-th gonality are all linearly equivalent to $l\cdot D$.
    \end{enumerate}
\end{conjecture*}

\begin{figure}[H]
\centering
\begin{tikzpicture}[scale=1.68, node/.style={circle, draw, fill=black, inner sep=1.2pt}]

% Graph 1

\def\y{0.1}
\node[node] (a1) at ({-5.5+\y},0) {};
\node[node] (a2) at ({-5.5+sqrt(3)/2+\y},0.5) {};
\node[node] (a3) at ({-5.5+\y},1) {};
\node[node] (a4) at ({-5.5-0.9396926207+\y},1.34202) {};
\node[node] (a5) at ({-5.5-0.9396926207+\y},-0.34202) {};

\draw (a1)--(a2)--(a3)--(a1);
\draw[bend left=18] (a1) to (a5);
\draw[bend right=18] (a1) to (a5);
\draw[bend left=18] (a3) to (a4);
\draw[bend right=18] (a3) to (a4);

% Graph 2

\node[node] (b1) at (-4,0) {};
\node[node] (b2) at (-3,0) {};
\node[node] (b3) at (-2,0) {};
\node[node] (b4) at (-1,0) {};
\node[node] (b5) at (-2,1) {};
 
\draw[bend left=18] (b1) to (b2);
\draw[bend right=18] (b1) to (b2);
\draw[bend left=18] (b3) to (b2);
\draw[bend right=18] (b3) to (b2);
\draw[bend left=24] (b3) to (b4);
\draw[bend right=24] (b3) to (b4);
\draw (b3) -- (b4);
\draw (b2) -- (b5);
\draw (b3) -- (b5);
\draw (b4) -- (b5);

% Graph 3

\def\x{-0.2}
\node[node] (v1) at ({0+\x},0) {};
\node[node] (v2) at ({1+\x},0) {};
\node[node] (v3) at ({0+\x},1) {};
\node[node] (v4) at ({1+\x},1) {};
\node[node] (v5) at ({1.9396926207+\x},1.34202) {};
\node[node] (v6) at ({0.539446+1.9396926207+\x},0.5) {};
\node[node] (v7) at ({1.9396926207+\x},-0.34202) {};
% Edges
\draw (v1) -- (v2);
\draw (v2) -- (v3);
\draw (v3) -- (v4);
\draw (v4) -- (v1);
\draw (v1) -- (v3);
\draw (v2) -- (v4);
\draw (v5) -- (v6);
\draw (v6) -- (v7);
\draw[bend left=30] (v5) to (v4);
\draw[bend right=10] (v5) to (v4);
\draw[bend left=10] (v5) to (v4);
\draw[bend right=30] (v5) to (v4);
\draw[bend left=30] (v7) to (v2);
\draw[bend right=10] (v7) to (v2);
\draw[bend left=10] (v7) to (v2);
\draw[bend right=30] (v7) to (v2);

\end{tikzpicture}
\caption{Examples of gonality-tight graphs satisfying the conjecture.}
\label{fig:examples}
\end{figure}

The conjecture is consistent with the known gonality sequence of the complete graphs computed in \cite{CoolsPanizzut_gonalityComplete}.
The first two parts of the conjecture appeared in \cite[Question 1.3]{FJKO_semigroupOfGonality}. In the same paper, it is  also shown that the first two parts hold for algebraic curves \cite[Lemma 7.2]{FJKO_semigroupOfGonality}. As we shall see, the third part holds as well for algebraic curves (\cref{lem:algebraicVersionOfConjecture}).

% \begin{lemma}\label{lem:algebraicVersionOfConjecture}
%     Let $C$ be a smooth curve with first gonality $k$ and second gonality $k+1$ for some integer $k$. Denote by $D$ a divisor that realises the second gonality. For $1\leq l\leq k-2$, a divisor realises the $\frac{l(l+3)}{2}$-the gonality if and only if it is equivalent to $l\cdot D$.
% \end{lemma}

% \begin{proof}
%    The result is an almost immediate consequence of \cite[Theorem 2.1]{Hartshorne_GorensteinAndNoether}, which Hartshorne attributes to Noether.  By \cite[Lemma 7.1]{FJKO_semigroupOfGonality}, we can assume that $C$ is a plane curve of degree $k+1$ and that the divisor obtained as a  line section of $C$ is equivalent to $D$.
   
%    Let $Z$ be a divisor realising the $\frac{l(l+3)}{2}$-th gonality. By \cite[Lemma 7.2]{FJKO_semigroupOfGonality}, the degree of $Z$ is $l\cdot (k+1)$. From Part (2.b) of Noether's theorem, the divisor $Z$ is obtained as the scheme-theoretic intersection of $C$ with another curve $C''$ of degree $l$. Therefore $Z$ is equivalent to $l\cdot D$. 
   
%    Similarly, since $D$ is a line section of $C$,  the divisor $l\cdot D$ is the intersection of $C$ with a curve of degree $l$, so the ``only if" direction of Noether's theorem implies that the rank of $l\cdot D$ is $\frac{l\cdot(l+3)}{2}$. 
% \end{proof}

\medskip

Next, we  ask for the following generalisation of \cref{mainTheorem:alternativeDefinitonForQBGraphs}.

\begin{question}
    Suppose that $G$ is a gonality-tight graph whose first gonality is equal to the number of vertices. Does it follow that $G$ is a quasi-banana graph?
\end{question}
\noindent We remark that the generalised banana graphs $B^*_{a,b}$ are \emph{not} gonality-tight when $b\neq a$, so they don't form a counter-example to the question. 

% \simun{I changed the ending of \cref{sec:introduction} to include the new result about subdividing edges. Let me know what you think and feel free to make any edits.} 
We conclude the paper by examining the result of subdividing the edges of quasi-banana graphs. 
%We show that equally subdividing edges  preserves the gonality-tight property, therefore giving a large family of gonality-tight graphs. 
Given a graph $G$ and an integer $\ell\geq1$, let $\sigma_\ell (G)$ denote the graph obtained from $G$ by subdividing every edge $\ell$ times. In general, $\gon_r (G)$ need not  equal $\gon_r(\sigma_ \ell (G))$ (see \cite[Theorem 1.4]{DBSW_gonalityCanBeDifferent}). However, in the case of quasi-banana graphs, the first two gonalities are preserved. 

\begin{maintheorem}\label{mainTheorem:subdividingEdges}
    Let $G$ be a quasi-banana graph with gonality $k$. For every integer $\ell\geq 1$, the graph $\sigma_\ell(G)$ is gonality-tight with gonality $k$. In particular, if $\Gamma$ is the metric graph obtained from $G$ by assigning length 1 to each edge, then $\Gamma$ is gonality-tight as well. 
    
    % \simun{Let $\Gamma(G)$ be the metric graph obtained by assigning length $1$ to every edge in a graph $G$. This is equivalent to saying that $\Gamma(G)$ is gonality-tight with gonality $k$. Is this a better way of phrasing the theorem?}
\end{maintheorem}

%It would also be interesting to explore the result of equally subdividing each of the edges. While we don't focus on this direction in the current paper, we expect (based on experimental evidence) this operation to preserve the gonality-tight property, giving an even larger family of gonality-tight graphs.

%\begin{question}\label{q:subdividingVertices}
    %Let $G$ be a quasi-banana graph and let $G'$ be the graph obtained by subdividing each edge of $G$. Is the graph  $G'$  gonality-tight as well? 
    
    % \simun{I think we can remove the part that every divisor on $G$ realising the second gonality continues to realise the second gonality on $G'$. I believe that it is true and follows from Theorem 3 in ``Rank of divisors on tropical curves" by Hladky, Kral, Norine. We can instead mention that the gonality sequence does not change.}\yoav{Good point.  I've also changed the statement into a question and the same with question 1.2.} \simun{I think an even more general result is true, the edges don't have to be subdivided an equal number of times, and the graph should stay gonality-tight. For example, we can subdivide one edge once, another edge twice and another five times, and not subdivide the remaining ones at all.}\yoav{Are there examples that support that?} \simun{The computational results support the more general version of the question. I also computed gonalities of some larger quasi-banana graphs with subdivided edges, and they remain gonality-tight.}
%\end{question}

It is now natural to ask whether this is true for gonality-tight graphs in general.

\begin{question}\label{question:subdividingEdges}
    Let $G$ be a gonality-tight graph with gonality $k$ and $\ell$ a positive integer. Is $\sigma_\ell(G)$ necessarily gonality-tight with gonality $k$?
\end{question}

Note that the graph obtained from equally subdividing the edges of a quasi-banana graph is \emph{no longer} a quasi-banana graph. 
Furthermore, subdividing the edges  \emph{unequally} may break gonality-tightness, as seen in \cref{example:unequal subdivision is not gonality tight}. 

\subsection{Strategy of the proof of the main theorems}

The backwards direction of \cref{mainTheorem:alternativeDefinitonForQBGraphs} is proven in \cref{sec:gonalityTightAndQBGaphs}. As a first step, we prove \cref{lemma:structureMainLemma}, which imposes conditions on the structure of gonality-tight graphs with a special divisor realising the second gonality. 
 In \cref{sec:quasiBananaGraphs}, we then repeatedly apply this lemma to show that all such graphs must be quasi-banana graphs. 
%The equivalence to quasi-banana graphs is then proven in \cref{sec:quasiBananaGraphs} by repeated use of the lemma. 
The proof of the forward direction is postponed to \cref{sec:ProofOfThmA}, as it relies on various tools developed in between. 
The proof of \cref{mainTheorem:propertiesOfQBGraphs}, which describes the properties of quasi-banana graphs, 
will be carried out  in \cref{sec:TwoTechnicalLemmas,sec:proofOfMainTheorem}: the former is dedicated to proving two lemmas that describe the effect on the gonalities of attaching a complete graph to a given graph and those lemmas will then be then used inductively to prove \cref{mainTheorem:propertiesOfQBGraphs} in \cref{sec:proofOfTheoremB}. 
Finally, we prove \cref{mainTheorem:subdividingEdges} in \cref{sec:proofOfTheoremC}.

\subsection*{Acknowledgements}
We thank the nameless referees for helpful comments on a previous version of this manuscript. 
We thank Dave Jensen for helpful discussions about the problem.
The project was undertaken as part of the Undergraduate Summer Research Programme at the School of Mathematics and Statistics in St Andrews and we thank them for their hospitality.
The research  was supported by the EPSRC New Investigator Award (grant number EP/X002004/1) and by the Student Research Bursaries of the Edinburgh Mathematical Society.

\section{Preliminaries} \label{sec: preliminaries}

We assume familiarity with basic terms in the theory of tropical divisors, such as degree, linear equivalence, $v$-reduced divisors, rank, and the Riemann--Roch theorem. 
There are many good sources such as \cite{BakerNorine2007, AC_rank, CLM_algebraicRank, LU_divisorsSurvey}, 
but our language mostly matches that of \cite[Section 2]{ADMYY_GonalitySequencesOfGraphs}. We do, however, briefly overview  some of the basic notions and establish some terminology.

Unless stated otherwise, all graphs are assumed to be finite connected loopless multigraphs. 
If $H$ is a subgraph of $G$, and $U\subseteq V(G)\setminus V(H)$, we write $H\cup U$ for the subgraph whose vertex set is $V(H) \cup U$ and whose edge set consists of all the edges in $H$ together with all edges that connect either two vertices in $U$ or a vertex in $U$ to a vertex in $V(H)$. Similarly, if $W\subseteq V(H)$, we write $H\setminus W$ for the graph obtained from $H$ by removing the vertices in $W$ and all edges incident to any vertex in $W$.

A \emph{divisor} on a graph $G$ is a $\mathbb{Z}$-linear sum of the vertices $V(G)$, that is, an element of the free abelian group generated by $V(G)$. If $H$ is a subgraph of $G$, we regard the group $\operatorname{Div}(H)$ of divisors on $H$ as a subgroup of $\operatorname{Div}(G)$, the group of divisors on $G$. In particular, if $u,v\in V(H)\subseteq V(G)$, we make no distinction between, e.g., $2u+v$ as a divisor on $H$, and the same element viewed as a divisor on $G$. A divisor is \emph{supported} on a subgraph $H$ if it has no chips outside of $H$, namely $D(u) = 0$ for all vertices $u$ not in $H$.

Since the rank of a divisor does depend on the ambient graph, we denote by $r_H(D)$ the rank of a divisor $D$  supported on $H$, thought of as a divisor on the subgraph $H$. 
If $D$ is a divisor on $G$ and $H$ is a subgraph on $G$, then we write $D|_H$ for the divisor $D$ restricted to the vertices of $H$. That is,
\[
D|_H=\sum_{w\in V(H)}D(w)\cdot w,
\]
where $D(w)$ denotes the coefficient of $w$ in $D$.  We similarly denote $f|_{H}$ for the restriction of a function $f$ to the vertices of $H$. 
We say that a divisor is \emph{equivalent to effective} or just \emph{winnable} if it is linearly equivalent to an effective divisor. 
Given a set of vertices $U$, we will use the phrase \emph{chip-firing from $U$} or \emph{set-firing $U$} to refer to chip-firing once from each vertex of $U$. 

We will make repeated use of the Dhar's burning algorithm. Given a divisor $D$, the algorithm can be described by starting a fire from a vertex $v$ that spreads under certain rules. If the fire doesn't burn the entire graph, there is an unburnt set of vertices $U$. We then chip-fire from $U$ and repeat the process. This process eventually terminates and the resulting divisor is the unique $v$-reduced divisor equivalent to $D$. See \cite[Section 2.3]{ADMYY_GonalitySequencesOfGraphs} and \cite[Section 2]{FJKO_semigroupOfGonality} for a detailed explanation.

The $k$-th \emph{gonality} of a graph, denoted $\gon_k(G)$, is the lowest degree of a divisor of rank $k$ for $k\geq 0$. We say that a divisor $D$ \emph{realises} the $k$-th gonality if it has rank $k$ and degree $\gon_k(G)$. A priori, there may be more than one divisor class realising any gonality. 
The \emph{gonality sequence} of a graph is the sequence $\gon_1(G),\gon_2(G),\ldots$. Note that the $0$-th  gonality is not interesting since it always equals $0$. Every gonality sequence is strictly increasing \cite[Lemma 3.1]{ADMYY_GonalitySequencesOfGraphs}. Of key interest to us are graphs whose second gonality exceeds the first by $1$, so we give them a name.
\begin{definition}\label{def:gonalityTight}
    Let $G$ be a graph with first gonality $k$ and second gonality $k+1$ for some $k$. Then $G$ is called a \emph{gonality-tight} graph.
\end{definition}

We denote the canonical divisor on a graph $H$ by $K_H$. The following is an immediate corollary of Riemann--Roch \cite{BakerNorine2007} that will be applied multiple times.

\begin{corollary}\label{corollary:corollariesOfRiemann-Roch}
    Let $H$ be a graph of genus $g$, and $D$ a divisor on $H$.  

    \begin{enumerate}
        \item If $\deg(D) = 2g - 2$, then 
        \[
        r_H(D) =
        \begin{cases}
            g - 1 & \text{if $D\sim K_H$,}\\
            g - 2 & \text{otherwise.}
        \end{cases}
        \]
        \item If $\deg(D)\geq 2g-1$, $r_H(D)=\deg(D)-g$.
        \item For $r\geq g$, $\gon_r(G)=g+r$.
    \end{enumerate}
\end{corollary}

\begin{proof}
    If $\deg(D) = 2g-2$ and $D\sim K_H$ then from Riemann--Roch, 
    \[
    r_H(D)=r_H(K_H)=r_H(K_H) - r_H(0) = \deg(K_H)-g+1 = g-1.
    \]
    If instead $D\not\sim K_H$, then 
    \[
    r_H(D) - r_H(K_H-D) = g-1.
    \]
    Since $D\not\sim K_H$, the divisor $K_H-D$ is of degree $0$, but not equivalent to the zero divisor. Its rank is therefore $-1$. It follows that $r_H(D) = g-2$.

    If $\deg(D)\geq 2g-1$, then 
    \[
    r_H(D) - r_H(K_H-D) = \deg(D) - g +1.
    \]
    Since the degree of $K_H-D$ is negative, it follows that $r_H(K_H-D) = -1$, so $r_H(D) = \deg(D)-g$.

    Finally, suppose that $r\geq g$. Since $r+g > 2g-1$, by Part (2), the rank of every divisor of degree $g+r$ equals $r$. From Part (2), the rank of any divisor of degree $2g-1$ is strictly smaller than $r$, so the gonality equals $r+g$. 
\end{proof}

\section{Gonality-tight and quasi-banana graphs}\label{sec:gonalityTightAndQBGaphs}

In this section, we explore  gonality-tight graphs whose second  gonality is realised by a special divisor. We show that such graphs always have a certain structure, and conclude that they coincide with quasi-banana graphs.

\subsection{Restricting the structure of gonality-tight graphs 
}\label{sec:gonalityTight}

Throughout, when  $v$ is a vertex of a graph $G$, and $G'$ is a subgraph of $G$, we denote by $N_{G'}(v)$ the set of neighbours of $v$ in $G'$.

\begin{lemma}\label{lemma:structureMainLemma}
    Let $G$ be a gonality-tight graph with a  vertex $v$, and suppose that the graph obtained after removing a vertex $v$ is a disjoint union of two  graphs $G'$ and $G''$ (where $G'$ and $G''$ are allowed to have multiple components   or be empty). 
    Suppose that, for some $k\geq 1$, the divisor 
    \begin{equation*}
        D_2 = (k+1)\cdot v + \sum_{w\in V(G')} w
    \end{equation*}
    is $v$-reduced and
    realises the second gonality.
    Then the following hold. 
    \begin{enumerate}
        \item There are exactly $k+1$ edges between $v$ and each of its neighbours in $G'$. 
             
        \item There is precisely one edge between every pair of neighbours of $v$ in $G'$. 
        
        \item If $V(G')\setminus N_{G'}(v)$ is non-empty, then there is a unique vertex $u\in N_{G'}(v)$ such that every path from $V(G')\setminus N_{G'}(v)$ to $v$  
        passes through $u$.

        \item 
        Set $n=|N_{G'}(v)|$. Then the $u$-reduced divisor equivalent to $D_2$ is
        \[
        (k+n+1)\cdot u + \sum_{w\in V(G')\setminus N_{G'}(v)}w.
        \]
    \end{enumerate}
\end{lemma}

The rest of this subsection is dedicated to  proving  the lemma.  
Throughout, we let $G$ be the graph with $\gon_1(G) + 1 = \gon_2(G)$, such that  one of the divisors realising the second gonality is a $v$-reduced divisor
\begin{equation*}
    D_2 = (k+1)\cdot v + \sum_{w\in V(G')} w,
\end{equation*}
with $k\geq 1$. This implies that $\gon_2(G) = \deg(D_2) = (k+1) + |V(G')|$ and $\gon_1(G) =\gon_2(G)-1= k+|V(G')|$. Note that,
\[
k=\gon_1(G)-|V(G')|\leq |V(G)|-|V(G')|=|V(G'')|+1.
\]
We denote by $u_1,\ldots,u_n$ the neighbours of $v$ in $G'$ (see  \cref{fig:1}). 
We start by proving two lemmas that will be used throughout.

\begin{figure}
\centering
\begin{tikzpicture}[scale=1, node/.style={circle, draw, fill=black, inner sep=1.2pt}]
    % v and graphs G', G''
    \node[node, label=above:{$v$}] (v) at (0,0) {};
    \draw (3.5,0) circle [radius=2];
    \node at (-3.5,0) {$G''$};
    \draw (-3.5,0) circle [radius=2];
    \node at (3.5,0) {$G'$};

    % Connecting v to G''
    \draw (v) to (-1.9,1.2);
    \draw (v) to (-1.9,-1.2);
    \draw (v) to (-1.55,0.44441);
    \draw[bend left=10] (v) to (-1.527,-0.32752);
    \draw[bend right=10] (v) to (-1.527,-0.32752);

    % Vertices u_i
    \node[node, label=above:{$u_1$}] (u1) at (2,1.32288) {};
    \node[node, label=right:{$u_2$}] (u2) at (1.617,0.67403) {};
    \node[node, label=right:{$u_3$}] (u3) at (1.5,0) {};
    \node[node, label=right:{$u_4$}] (u4) at (1.617,-0.67403) {};
    \node[node, label=below:{$u_5$}] (u5) at (2,-1.32288) {};

    % Connecting v to G'
    \draw[bend right=5] (v) to (u1);
    \draw[bend left=5] (v) to (u1);
    \draw (v) to (u2);
    \draw (v) to (u3);
    \draw[bend right=5] (v) to (u4);
    \draw[bend left=5] (v) to (u4);
    \draw (v) to (u5);
    
\end{tikzpicture}
\caption{An illustration of the graph $G$ from \cref{lemma:structureMainLemma}.}
\label{fig:1}
\end{figure}

\begin{lemma}\label{lemma:canRestrictToComponent}
    Let $H$ be a graph  with a cut vertex $v$ whose removal results in the disjoint union of (not necessarily connected) subgraphs $H'$ and $H''$. Suppose that $D$ is a divisor of rank $r$ with the property that $ D|_{H'} - v$ has rank $-1$ (as a divisor on $H'\cup\{v\}$). 
Then $D|_{H''\cup \{v\}}$ has rank at least $r$ as a divisor on $H''\cup\{v\}$. 
\end{lemma}

\begin{proof}
The lemma  follows as a special case of \cite[Remark 4.9]{AminiBaker_metrizedComplexes}, however, we provide an independent proof as follows. 
    Let $E$ be an effective divisor of degree $r$ on $H''\cup\{v\}$. 
    As $r_H(D) = r$, let $f$ be the piecewise linear function on $H$ such that $D-E+\ddiv(f)\geq 0$. We will show that the divisor $D|_{H''\cup\{v\}} - E +\ddiv(f|_{H''\cup\{v\}})$ is effective.  This divisor can only differ from $D-E+\ddiv(f)$ at $v$, so we only need to check that it is effective at $v$. The difference between $\ddiv(f)(v)$ and $ \ddiv(f|_{H''\cup\{v\}})(v)$ is exactly the sum  slopes of $f$ at $v$ incoming from $H'$. 
    Therefore, if the divisor $D|_{H''\cup\{v\}} - E +\ddiv(f|_{H''\cup\{v\}})$ is not effective at $v$, 
    the sum of the slopes of $f$ at $v$ incoming from $H'$ must be positive. It then follows that $\ddiv(f|_{H'\cup \{v\}})$ is positive at $v$, so $-v + D|_{H'} + \ddiv(f|_{H'\cup\{v\}})$ is effective, which is a contradiction. 
\end{proof}

% The lemma above also follows as a special case of \cite[Remark 4.9]{AminiBaker_metrizedComplexes}
%by choosing $\Gamma' = H_2\cup\{v\}$ and $\Gamma_2 = H'\cup\{v\}$ (where $\Gamma_1$ and $\Gamma_2$ are in the terminology of \cite{AminiBaker_metrizedComplexes}).

\begin{lemma}\label{lemma:restrictedRank}
    Let $H$ be a graph  with a cut vertex $v$ whose removal results in the disjoint union of (not necessarily connected) subgraphs $H'$ and $H''$. Let $D$ be a divisor on $H''\cup\{v\}$ with rank $r$. Let $F$ be an effective divisor on $H$, such that $D\leq F|_{H''\cup\{v\}}$ 
    Then, whenever $E$ is an effective divisor with degree $r$ supported on $H''\cup\{v\}$, the divisor $F-E$ is equivalent to an effective divisor on $H$.
\end{lemma}
\begin{proof}
  
    Since $r_{H''\cup \{v\}}(D) = r$, there is a piecewise linear function $f''$ on $H''\cup\{v\}$ such that $D - E +\ddiv(f'')\geq 0$.  Let $f$ be the piecewise linear function on $G$ that is constant on $H'\cup \{v\}$ and coincides with $f''$ on $H''\cup\{v\}$.  Since $D\leq F|_{H''\cup\{v\}}$ and $f$ is constant on $H'\cup\{v\}$, it follows that $F-E + \ddiv(f)\geq 0$.
\end{proof}

We now return to our fixed graph $G$ and divisor $D_2$.

\begin{claim}\label{claim:kChipsOnv}
    The divisor $k\cdot v$ has rank at least $1$ as a divisor on $G''\cup\{v\}$. In particular, if $D$ is an effective divisor  on $G$ such that $D(v)\geq k$,  the divisor $D-w$ is equivalent to an effective divisor for every vertex  $w\in V(G'')$. 
    
\end{claim}
\begin{proof}
    
    Since our fixed divisor $D_2$ realises the second gonality, the divisor $D_2-v$ has rank $1$.
    Moreover, $D_2$ is $v$-reduced, and so the divisor 
    \[
    D_2|_{G'}-v=-v+\sum_{w\in V(G')} w = (D_2- v)|_{G'}-v
    \]
    is as well. In particular, it has rank $-1$ as a divisor on $G'\cup\{v\}$. Applying \cref{lemma:canRestrictToComponent} to $D_2-v$,
    we conclude that $k\cdot v=(D_2-v)|_{G''\cup\{v\}}$ has rank at least $1$ as a divisor on $G''\cup\{v\}$. Now, by \cref{lemma:restrictedRank}, for every $w\in V(G'')$, the divisor $k\cdot v-w$ is equivalent to an effective divisor on $G$. As $D-w\geq k\cdot v-w$, the claim is proven.
\end{proof}

Recall that $u_1,u_2,\ldots,u_n$ are the neighbours of $v$ in $G'$.

\begin{figure}
\centering
\begin{tikzpicture}[scale=0.8, node/.style={circle, draw, fill=black, inner sep=1.2pt}]
    % v and u_i, and edges between them
    \node[node] (v) at (0,2.6) {};
    \node[above left=-0.15cm and 0.02cm of v] {$v$};
    \node[node, fill=red, draw=red, label=below:{$u_i$}] (ui) at (0,0) {};
    \draw[bend left=18] (v) to (ui);
    \draw[bend left=6] (v) to (ui);
    \draw[bend right=6] (v) to (ui);
    \draw[bend right=18] (v) to (ui);
    \node at (-0.5,1.3) {$\ell$};
    
    % Subgraph G'' and connecting v to G''
    \draw (0,4) circle [radius=0.7];
    \node at (0,4) {$G''$};
    \draw (v) -- (-0.61, 3.65663);
    \draw (v) -- (0.61, 3.65663);
    \draw (v) -- (-0.19, 3.32628);
    \draw (v) -- (0.19, 3.32628);

    % Burnt and unburnt parts of G'
    \draw (-4,0.7) ellipse [x radius=1.8, y radius=3];
    \node at (-4,0.7) {$U$};
    \draw[red, thick] (4,0.7) ellipse [x radius=1.8, y radius=3];
    \node at (4,0.7) {$B$};

    % Connecting u_i to B
    \draw[draw=red] (ui) -- (2.28, -0.18443);
    \draw[draw=red] (ui) -- (2.36, -0.53648);
    \draw[draw=red, bend left=8] (ui) to (2.224, 0.21174);
    \draw[draw=red, bend right=8] (ui) to (2.224, 0.21174);

    % Vertices w_i and connecting them to v
    \node[node, fill=red, draw=red, label=below right:{$w_3$}] (w3) at (2.26,1.46811) {};
    \node[node, fill=red, draw=red] (w2) at (2.347,1.88742) {};
    \node[above right=-0.4cm and 0.02cm of w2] {$w_2$};
    \node[node, fill=red, draw=red] (w1) at (2.48,2.30693) {};
    \node[above right=-0.25cm and 0.02cm of w1] {$w_1$};
    \draw[bend left=7] (v) to (w1);
    \draw[bend right=7] (v) to (w1);
    \draw (v) to (w3);
    \draw (v) to (w2);
    \node at (1.35,2.83) {$\ell'$};

    % Edges between U and B
    \draw (-2.444,-0.80819) -- (2.53, -1.03133);
    \draw (-2.57,-1.12201) -- (2.53, -1.03133);
    \draw (-2.735,-1.43423) -- (2.606, -1.19793);
    
    % Connecting v to U
    \draw (v) to (-2.43, 2.16733);
    \draw (v) to (-2.306, 1.71429);

\end{tikzpicture}
\caption{The result of starting a fire at $u_i$ with respect to $D_2-2u_i$. Parts of the graph that burn down  are highlighted in red.}
\label{fig:2}
\end{figure}
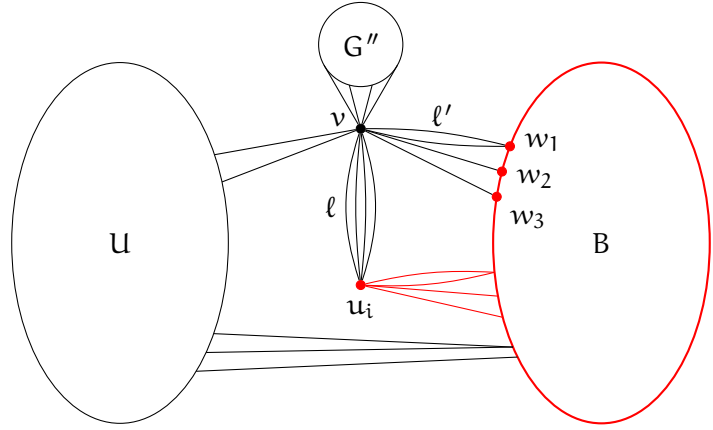

\begin{claim}\label{claim:edgesBetweenvAndui}
    There are exactly $k+1$ edges between $v$ and each $u_i$.
\end{claim}
\begin{proof}
    Fix one of the vertices $u_i$.  
    As $D_2$ is $v$-reduced, if $v$ is ever on fire, the whole graph burns down. Since $D_2$ has rank $2$, the divisor $D_2-2\cdot u_i$ is winnable. In particular, if we start a fire from  $u_i$, it cannot burn $v$, so there must be  fewer than $k+2$ edges between $u_i$ and $v$. Denote by  $\ell\leq k+1$  the number of edges between $v$ and $u_i$.

    Now, start a fire from $u_i$ with respect to the divisor $D_2$. Denote by $B$ the subgraph of  $ G'\setminus\{u_i\}$ that is burned by the fire and by $U$ the subgraph that isn't. The subgraphs $B$ and $U$ may have multiple components. Let $\ell'$ be the number of edges between $v$ and the vertices on fire other than $u_i$ (see \cref{fig:2}). Denote by $w_1,w_2,\ldots, w_m$ the vertices at the other end of those edges. Note that the number of those vertices might be smaller than $\ell'$ due to multiple edges.

We claim that $\ell'=0$. Otherwise, we will show that the divisor $D_2 - u_i - \sum w_j$ has rank $1$, which is a contradiction since $\deg(D_2 - u_i - \sum w_j)<\deg(D_2)-1=\gon_1(G)$. Indeed, if $u\in V(G')\cup\{v\}$ is any vertex other than $u_i$ or $w_1, w_2,\ldots, w_m$, then $D_2 - u_i - \sum w_j$ already has a chip at $u$. If $u\in V(G'')$ then, by \cref{claim:kChipsOnv}, the divisor $D_2 - u_i - \sum w_j-u$ is equivalent to an effective divisor on $G$. 
Now, suppose that $u$ is  either $u_i$ or one of $w_1,w_2\ldots w_m$. The fact that $U\cup\{v\}$ didn't burn by the fire starting from $u_i$ precisely means that chip-firing from this set doesn't create any new debts. But chip-firing from the set also brings a chip to $u$.

It follows that $\ell'=0$ and  the number of edges between $B\cup \{u_i\}$ and $v$ is $\ell$. We will show that, if $\ell<k+1$ then the divisor $D_2 - v - u_i$ has rank $1$, which is a contradiction. So suppose that $\ell< k+1$. The divisor $D_2-v-u_i$ already has a chip on  on every vertex of $G'\cup\{v\}$ other than $u_i$. If $w\in V(G'')$ then by \cref{claim:kChipsOnv}, the divisor $D_2 - v - u_i-w$ is equivalent to an effective divisor. We are left with checking that the divisor obtained by taking a chip from $u_i$ is winnable. Start a fire from $u_i$. Since the number of edges between $B\cup \{u_i\}$ and $v$ is $\ell< k+1$, the fire doesn't burn $v$. Chip-firing the unburnt part results in clearing the debt at $u_i$. We have reached a contradiction, and the only possibility is that  $k+1\leq \ell < k+2$, so the claim is proven.
\end{proof}

\begin{figure}
\centering
\begin{tikzpicture}[scale=3.5, node/.style={circle, draw, fill=black, inner sep=1.2pt}]
    % Vertex v
    \node[node, label=above:{$v$}] (v) at (0.25,0) {};

    % Subgraph G''
    \draw (-0.45,0) circle [radius=0.45];
    \node at (-0.45,0) {$G''$};

    % Edges connecting v to G''
    \draw[bend left=10] (v) to (-0.01,0.09434);
    \draw[bend right=10] (v) to (-0.01,0.09434);
    \draw (v) to (-0.01,-0.09434);
    \draw (v) to (-0.16,0.34409);         
    \draw (v) to (-0.16, -0.34409);       

    % Part of G' that is not on fire
    \draw (1.1,0) ellipse [x radius=0.47, y radius=0.65];

    % u_i-s
    \node[node, label=below right:{$u_1$}] (u1) at (0.73,0.40083) {};
    \node[node, label=right:{$u_2$}] (u2) at (0.62,0) {};
    \node[node, label=above right:{$u_3$}] (u3) at (0.73,-0.40083) {};

    % Edges between v and u_i-s
    \draw[bend left=10] (v) to (u1);
    \draw (v) -- (u1);
    \draw[bend right=10] (v) to (u1);
    \draw[bend left=10] (v) to (u2);
    \draw (v) -- (u2);
    \draw[bend right=10] (v) to (u2);
    \draw[bend left=10] (v) to (u3);
    \draw (v) -- (u3);
    \draw[bend right=10] (v) to (u3);

    % Part of the graph on fire
    \draw[red, thick] (2.4,0) ellipse [x radius=0.47, y radius=0.65];

    % Nodes w_i-s that are on fire
    \node[node, fill=red, draw=red, label=above left:{$w_1$}] (w1) at (2.05,0.43382) {};
    \node[node, fill=red, draw=red, label=right:{$w_2$}] (w2) at (1.94,0.13337) {};
    \node[node, fill=red, draw=red, label=right:{$w_3$}] (w3) at (1.94,-0.13337) {};
    \node[node, fill=red, draw=red, label=below left:{$w_4$}] (w4) at (2.05,-0.43382) {};

    % Connecting w_i-s to the unburnt part
    \draw (w1) -- (1.42,0.47607);
    \draw (w2) -- (1.52,0.29174);
    \draw (w2) -- (1.57,0);
    \draw (w3) -- (1.55,-0.1876);
    \draw (w4) -- (1.43,-0.46283);

    % vertex w
    \node[node, fill=red, draw=red, label=above:{$w$}] (w) at (2.5,0) {};
\end{tikzpicture}
\caption{Graph $G$ after a fire is started at $w$ with respect to $D_2$. Parts of the graph on fire are highlighted in red.}
\label{fig:3}
\end{figure}
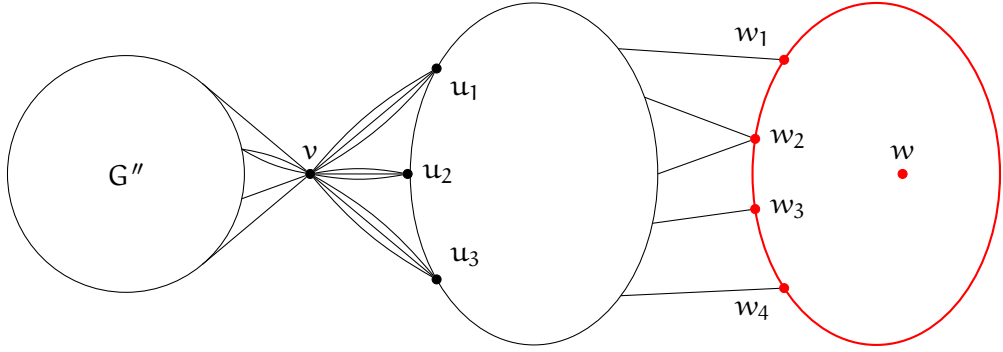

\begin{claim}\label{claim:exactlyOneu_iBurns}
    If a fire is started at some $w\in V(G')$ with respect to $D_2$, then exactly one $u_i\in N_{G'}(v)$ burns down.
\end{claim}
\begin{proof}

If $w\in N_{G'}(v)$, then it follows from the proof of \cref{claim:edgesBetweenvAndui} that starting a fire at $w$ does not cause any of the other neighbours of $v$ to catch fire. Hence the claim holds in this case. 
    
Now, let $w\in V(G')\setminus N_{G'}(v)$. Suppose, for contradiction, that no neighbour of $v$ burns down when a fire is started at $w$ with respect to $D_2$, as illustrated in \cref{fig:3}. Let $w_1,w_2,\ldots, w_m$ be the vertices that burn, but have at least one unburnt neighbour. 
    
We claim that the divisor $D_2-v-\sum w_i$ has rank at least $1$, which is a contradiction since its degree is smaller than the gonality.
If $w'\in V(G'')$, then by \cref{claim:kChipsOnv}, the divisor $D_2-v-\sum w_i-w'$ is equivalent to an effective divisor. If $w'$ is one of the $w_1,w_2,\ldots, w_m$, then set-firing the unburnt part of $G$ in \cref{fig:3} adds at least one chip to $w'$, without creating any debt elsewhere. $D_2-v-\sum w_i-w'$ is again equivalent to effective, and so the rank of $D_2-v-\sum w_i$ is at least $1$. 

Therefore, there must be some $u_i\in N_{G'}(v)$ that burns  when a fire is started at $w$. By \cref{claim:edgesBetweenvAndui}, there are $k+1$ edges between $v$ and each of its neighbours in $G'$, so if another $u_j\neq u_i$ burns down, the fire will spread to $v$ and then the whole of $G$. This would imply that $D_2-2\cdot w$ is $w$-reduced and therefore unwinnable, contradicting that $D_2$ has rank $2$. Hence, exactly one vertex $u_i$ burns, as claimed.
\end{proof}

We say that a vertex $u$ \emph{separates} $w$ from $v$ if every path from $v$ to $w$ passes through $u$. As we are mostly interested in separating vertices from the fixed vertex $v$, we will often just say that $u$ separates $w$, without mentioning $v$.
For our fixed graph $G$, we will refer to the set of vertices separated by one of the vertices $u_i$ 
as the \emph{branch} of $u_i$, and denote it by $B_i$.

\begin{claim}\label{claim:branches}
    Setting fire with respect to $D_2$ to $u_i$, or any vertex in its branch $B_i$,  results in burning precisely $\{u_i\}\cup B_i$.
\end{claim}

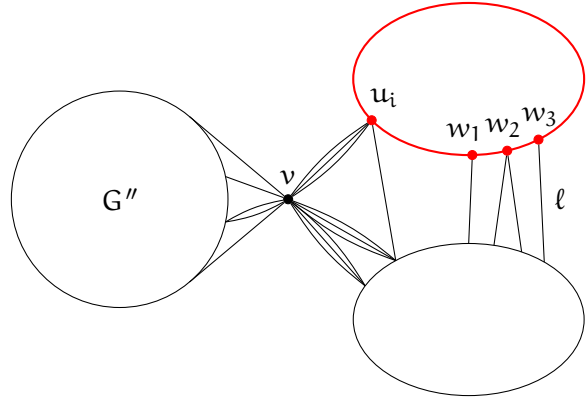
\begin{figure}
\centering
\begin{tikzpicture}[scale=3.2, node/.style={circle, draw, fill=black, inner sep=1.2pt}]
    % Vertex v
    \node[node, label=above:{$v$}] (v) at (0.25,0) {};

    % Subgraph G''
    \draw (-0.45,0) circle [radius=0.45];
    \node at (-0.45,0) {$G''$};

    % Edges connecting v to G''
    \draw[bend left=10] (v) to (-0.01,-0.09434);
    \draw[bend right=10] (v) to (-0.01,-0.09434);
    \draw (v) to (-0.01,0.09434);
    \draw (v) to (-0.16,0.34409);         
    \draw (v) to (-0.16, -0.34409);

    % Part of G' on fire
    \draw[red, thick] (1,0.5) ellipse [x radius=sqrt(0.23), y radius=sqrt(0.1)];
    \node[node, fill=red, draw=red] (ui) at (0.597, 0.32857) {};
    \node[above right=0cm and -0.2cm of ui] {$u_i$};
    \draw (v) to (ui);
    \draw[bend left=10] (v) to (ui);
    \draw[bend right=10] (v) to (ui);
    \node[node, fill=red, draw=red] (w1) at (1.015, 0.183927) {};
    \node[above right=0cm and -0.52cm of w1] {$w_1$};
    \node[node, fill=red, draw=red] (w2) at (1.16, 0.20189) {};
    \node[above right=0cm and -0.47cm of w2] {$w_2$};
    \node[node, fill=red, draw=red] (w3) at (1.29, 0.24814) {};
    \node[above right=0cm and -0.4cm of w3] {$w_3$};
    
    % Part of G' not on fire
    \draw (1,-0.5) ellipse [x radius=sqrt(0.23), y radius=sqrt(0.1)];
    \draw (v) to (0.697, -0.25488);
    \draw[bend left=10] (v) to (0.697, -0.25488);
    \draw[bend right=10] (v) to (0.697, -0.25488);
    \draw (v) to (0.57,-0.35997);
    \draw[bend left=10] (v) to (0.57,-0.35997);
    \draw[bend right=10] (v) to (0.57,-0.35997);

    % Connecting parts of G'
    \draw (ui) to (0.697, -0.25488);
    \draw (w1) to (1,-0.183772);
    \draw (w2) to (1.22, -0.219008);
    \draw (w2) to (1.105, -0.191444);
    \draw (w3) to (1.314,-0.260977);
    \node at (1.38,0) {$\ell$};

\end{tikzpicture}
\caption{The graph $G$ after starting a fire at $u_i$ with respect to $D_2$. Here, $\ell$ is the number of edges between vertices on fire other than $u_i$ and vertices of $G$ not on fire. Burnt regions are highlighted in red.}
\label{fig:4}
\end{figure}

\begin{proof}
Start a fire at $u_i$ with respect to $D_2$.
For the rest of the proof, we will use the notations from \cref{fig:4}. Our goal is to show that the vertices that burn are exactly the ones separated by $u_i$, which is equivalent to showing that $\ell=0$. So assume that $\ell > 0$. Let $w_1,w_2, \ldots w_m$ be the vertices of the burned locus that are adjacent to the unburnt locus.  Similarly to previous arguments, we will reach a contradiction by showing that the divisor $D_2 - u_i - \sum w_j$ has rank $1$. 
If we remove a chip from $w\in V(G'')$, then by \cref{claim:kChipsOnv}, $D_2 - u_i - \sum w_j-w$ is equivalent to an effective divisor. It suffices to check that, after removing a chip from $u_i$ or one of the vertices $w_j$, the divisor is equivalent to effective. 
So take a chip away from one of those vertices, which we denote $w$, and start a fire from $u_i$. Since the vertices $u_i,w_1,w_2,\ldots, w_m$ where burnt when starting a fire from $u_i$ with respect to the divisor $D_2$, the burnt locus will not change when starting a fire from $u_i$ with respect to $D_2 - u_i - \sum w_j - w$. Now, chip-firing from the complement of the burnt locus results in an effective divisor.

Next, let $w$ be a vertex in $B_i$. By \cref{claim:exactlyOneu_iBurns}, after starting a fire at $w$ with respect to $D_2$, exactly one neighbour of $v$ will burn down. As $u_i$ separates $w$ from all other neighbours of $v$, that is precisely the neighbour that catches fire. Once $u_i$ burns, the fire spreads through all of $\{u_i\}\cup B_i$ and no further, as desired.
\end{proof}

\begin{claim}\label{claim:oneEdgeBetweenu_iAndu_j}
    There is exactly one edge between $u_i$ and $u_j$ for any $i$ and $j$. 
\end{claim}

\begin{proof}
    If there were two or more edges between $u_i$ and $u_j$, then setting fire from $u_i$ with respect to $D_2$ would burn $u_j$, contradicting \cref{claim:exactlyOneu_iBurns}.
    It follows that  there is at most one edge between each  $u_i$ and $u_j$.

    Now, assume there are vertices $u_i$ and $u_j$ with no edges between them. Consider the divisor $D_2 - u_i - u_j$. If we take a chip away from $u_i$ and set it on fire, since there is no edge between $u_i$ and $u_j$, by \cref{claim:branches}, the only part of the graph other than $u_i$ that will burn down is $B_i$. By chip-firing the unburnt vertices, including $v$, we obtain an effective divisor. 
    Similarly, taking away a chip from $u_j$ does not make the divisor unwinnable.
Finally, for any $w\in V(G'')$, $D_2 - u_i - u_j-w$ is equivalent to an effective divisor, by \cref{claim:kChipsOnv}. 
Therefore, $D_2 - u_i - u_j$ has rank at least $1$, which is a contradiction as its degree is smaller than the gonality of $G$. 
\end{proof}

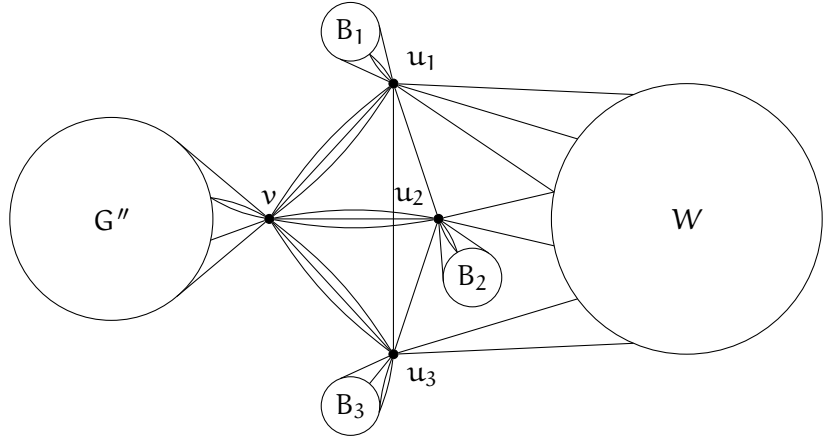
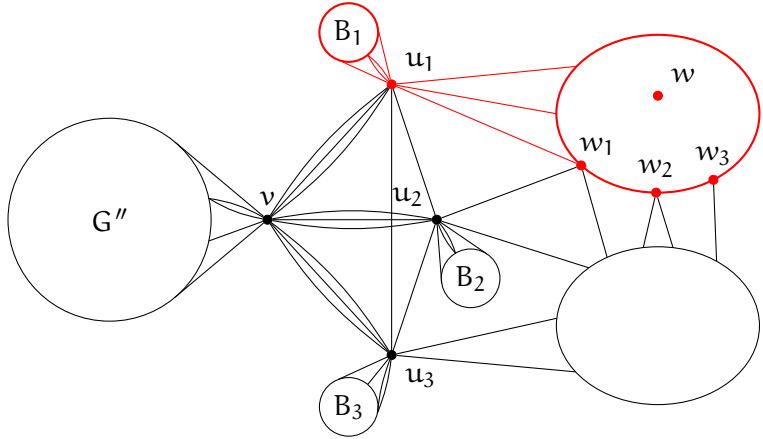
\begin{figure}
\centering
\begin{subfigure}[b]{0.9\textwidth}
\centering
\begin{tikzpicture}[scale=3, node/.style={circle, draw, fill=black, inner sep=1.2pt}]
    % Vertex v
    \node[node, label=above:{$v$}] (v) at (0.25,0) {};

    % Subgraph G''
    \draw (-0.45,0) circle [radius=0.45];
    \node at (-0.45,0) {$G''$};

    % Edges connecting v to G''
    \draw[bend left=10] (v) to (-0.01,0.09434);
    \draw[bend right=10] (v) to (-0.01,0.09434);
    \draw (v) to (-0.01,-0.09434);
    \draw (v) to (-0.16,0.34409);         
    \draw (v) to (-0.16, -0.34409);

    % u_is and edges between v and u_is, and between u_is
    \node[node, label=above left:{$u_2$}] (u2) at (1,0) {};
    \node[node, label=above right:{$u_1$}] (u1) at (0.8,0.6) {};
    \node[node, label=below right:{$u_3$}] (u3) at (0.8,-0.6) {};
    \draw[bend left=10] (v) to (u1);
    \draw (v) -- (u1);
    \draw[bend right=10] (v) to (u1);
    \draw[bend left=10] (v) to (u2);
    \draw (v) -- (u2);
    \draw[bend right=10] (v) to (u2);
    \draw[bend left=10] (v) to (u3);
    \draw (v) -- (u3);
    \draw[bend right=10] (v) to (u3);
    \draw (u1) -- (u2);
    \draw (u1) -- (u3);
    \draw (u2) -- (u3);

    % Drawing B_1
    \draw (0.61,0.83) circle [radius=0.13];
    \node at (0.61,0.83) {$B_1$};
    \draw (u1) to (0.737, 0.857767);
    \draw (u1) to (0.564, 0.708411);
    \draw[bend right=10] (u1) to (0.694, 0.730783);
    \draw[bend left=10] (u1) to (0.694, 0.730783);
    
    % Drawing B_2
    \draw (1.15,-0.26) circle [radius=0.13];
    \node at (1.15,-0.26) {$B_2$};
    \draw (u2) to (1.022, -0.282716);
    \draw (u2) to (1.222, -0.15176);
    \draw[bend left=10] (u2) to (1.083, -0.148595);
    \draw[bend right=10] (u2) to (1.083, -0.148595);
    
    % Drawing B_3
    \draw (0.61,-0.83) circle [radius=0.13];
    \node at (0.61,-0.83) {$B_3$};
    \draw[bend right=7] (u3) to (0.737, -0.857767);
    \draw[bend left=7] (u3) to (0.737, -0.857767);
    \draw (u3) to (0.564, -0.708411);
    \draw (u3) to (0.694, -0.730783);
    
    % Drawing W
    \draw (2.1,0) circle [radius=0.6];
    \node at (2.1,0) {$W$};
    \draw (u2) to (1.512,-0.119398);
    \draw (u2) to (1.512,0.119398);
    \draw (u1) to (1.512,0.119398);
    \draw (u1) to (1.864, 0.551638);
    \draw (u1) to (1.616, 0.354604);
    \draw (u3) to (1.616, -0.3546);
    \draw (u3) to (1.864, -0.55164);
\end{tikzpicture}
\caption*{(a) Sketch of the graph $G$ including $W$ and branches of each $u_i$.}
\end{subfigure}

\vspace{1em} % Space between subfigures

\begin{subfigure}[b]{0.9\textwidth}
\centering
\begin{tikzpicture}[scale=3, node/.style={circle, draw, fill=black, inner sep=1.2pt}]
    % Vertex v
    \node[node, label=above:{$v$}] (v) at (0.25,0) {};

    % Subgraph G''
    \draw (-0.45,0) circle [radius=0.45];
    \node at (-0.45,0) {$G''$};

    % Edges connecting v to G''
    \draw[bend left=10] (v) to (-0.01,0.09434);
    \draw[bend right=10] (v) to (-0.01,0.09434);
    \draw (v) to (-0.01,-0.09434);
    \draw (v) to (-0.16,0.34409);         
    \draw (v) to (-0.16, -0.34409);

    % u_is and edges between v and u_is, and between u_is
    \node[node, label=above left:{$u_2$}] (u2) at (1,0) {};
    \node[node, fill=red, draw=red, label=above right:{$u_1$}] (u1) at (0.8,0.6) {};
    \node[node, label=below right:{$u_3$}] (u3) at (0.8,-0.6) {};
    \draw[bend left=10] (v) to (u1);
    \draw (v) -- (u1);
    \draw[bend right=10] (v) to (u1);
    \draw[bend left=10] (v) to (u2);
    \draw (v) -- (u2);
    \draw[bend right=10] (v) to (u2);
    \draw[bend left=10] (v) to (u3);
    \draw (v) -- (u3);
    \draw[bend right=10] (v) to (u3);
    \draw (u1) -- (u2);
    \draw (u1) -- (u3);
    \draw (u2) -- (u3);

    % Drawing B_1
    \draw[red, thick] (0.61,0.83) circle [radius=0.13];
    \node at (0.61,0.83) {$B_1$};
    \draw[draw=red] (u1) to (0.737, 0.857767);
    \draw[draw=red] (u1) to (0.564, 0.708411);
    \draw[draw=red, bend right=10] (u1) to (0.694, 0.730783);
    \draw[draw=red, bend left=10] (u1) to (0.694, 0.730783);
    
    % Drawing B_2
    \draw (1.15,-0.26) circle [radius=0.13];
    \node at (1.15,-0.26) {$B_2$};
    \draw (u2) to (1.022, -0.282716);
    \draw (u2) to (1.222, -0.15176);
    \draw[bend left=10] (u2) to (1.083, -0.148595);
    \draw[bend right=10] (u2) to (1.083, -0.148595);
    
    % Drawing B_3
    \draw (0.61,-0.83) circle [radius=0.13];
    \node at (0.61,-0.83) {$B_3$};
    \draw[bend right=7] (u3) to (0.737, -0.857767);
    \draw[bend left=7] (u3) to (0.737, -0.857767);
    \draw (u3) to (0.564, -0.708411);
    \draw (u3) to (0.694, -0.730783);
    
    % Drawing burnt part of W
    \draw[red, thick] (1.98,0.47) ellipse [x radius=0.45, y radius=0.35];
    \node[node, fill=red, draw=red, label=above right:{$w$}] (w) at (1.98,0.55) {};
    
    \node[node, fill=red, draw=red] (w1) at (1.64, 0.24072) {};
    \node[above right=-0.05cm and -0.2cm of w1] {$w_1$};
    %\node[node, fill=red, draw=red] (w1) at (2.347,1.88742) {};
    %\node[above right=-0.4cm and 0.02cm of w2] {$w_2$};
    
    \node[node, fill=red, draw=red, label=above:{$w_2$}] (w2) at (1.972, 0.12006) {};
    \node[node, fill=red, draw=red, label=above:{$w_3$}] (w3) at (2.226, 0.17693) {};
    
    % Connecting u_is to W
    \draw[draw=red] (u1) to (1.53, 0.47);
    \draw[draw=red] (u1) to (1.62, 0.68);
    \draw[draw=red] (u1) to (w1);
    \draw (u2) to (w1);
    \draw (u2) to (1.674,-0.21338);
    
    % Edges between unburnt and burnt parts of W
    \draw (w1) to (1.755, -0.16689);
    \draw (w2) to (1.915, -0.12367);
    \draw (w2) to (2.047, -0.1239);
    \draw (w3) to (2.24, -0.18433);
    
    % Drawing unburnt part of W
    \draw (1.98,-0.47) ellipse [x radius=0.45, y radius=0.35];
    \draw (u3) to (1.615,-0.67471);
    \draw (u3) to (1.532, -0.43704);
\end{tikzpicture}
\caption*{(b) Graph $G$ after a fire is started at $w$ with respect to $D_2$.} 

\end{subfigure}

\caption{Graph $G$ before and after the vertex $w$ is set on fire. The parts of $G$ that burn down are highlighted in red.}
\label{fig:5}
\end{figure}

\begin{claim}\label{claim:WIsEmpty}
   Let $W$ be the complement in $V(G')$ of the union of the sets $\{u_i\}\cup B_i$ for $i=1,\ldots, n$. That is , $W$ is  the vertices in $G'$ that don't burn when starting a fire from any $u_i$ with respect to $D_2$ (See \cref{fig:5}(a) for an illustration).
   Then $W$ is empty.
\end{claim}

\begin{proof}
    Assume for a contradiction that $W$ is not empty, so it contains at least one vertex $w$. 
    Start a fire at $w$ with respect to $D_2$.
    By \cref{claim:exactlyOneu_iBurns}, the fire will burn exactly one neighbour $u_i$ of $v$ and therefore doesn't reach $v$, as depicted in \cref{fig:5}(b).
    By assumption, $w$ is not in a branch of $u_i$, so there exist some $w_1,w_2\ldots,w_m\in W$ such that setting a fire from $w$ burns them, and they are adjacent to some unburnt vertex.

    Now, consider the divisor $D_2 - u_i - \sum w_j$. We claim that its rank is at least  $1$, which is a contradiction, since its degree is lower than the gonality of $G$. If a chip is taken from some \(w \in V(G'')\), then by \cref{claim:kChipsOnv}, the divisor $D_2 - u_i - \sum w_j - w$ is still winnable.
We now only need to check that the divisor is equivalent to effective after taking away a chip from $u_i$ or one of the vertices $w_j$.  So take away a chip from one of those vertices, that we denote by $w'$, and start a fire from it. 
    Since we have only taken away chips from vertices that burnt when starting a fire from $u_i$ with respect to $D_2$, the set of vertices that catch fire is going to be a subset of the ones that are on fire in \cref{fig:5}(b). In particular, the vertex $w'$  will have an unburnt neighbour. Therefore, chip-firing from  the set of unburnt vertices will remove the debt from $G$.
\end{proof}

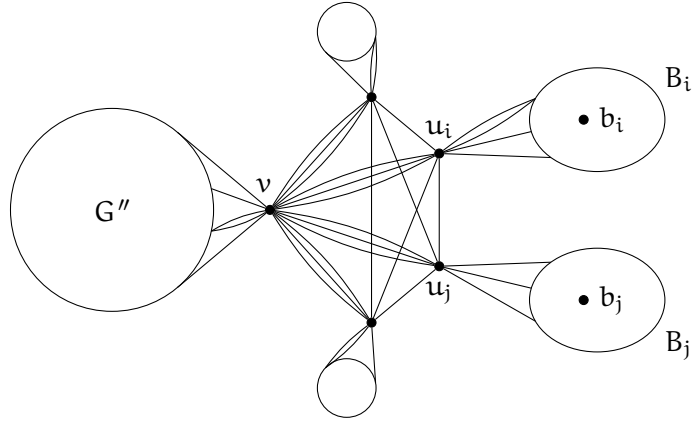
\begin{figure}
\centering
\begin{tikzpicture}[scale=3, node/.style={circle, draw, fill=black, inner sep=1.2pt}]
    % Vertex v
    \node[node] (v) at (0.25,0) {};
    \node[above left=0.06cm and -0.2cm of v] {$v$};
    
    % Subgraph G''
    \draw (-0.45,0) circle [radius=0.45];
    \node at (-0.45,0) {$G''$};

    % Edges connecting v to G''
    \draw[bend left=10] (v) to (-0.01,-0.09434);
    \draw[bend right=10] (v) to (-0.01,-0.09434);
    \draw (v) to (-0.01,0.09434);
    \draw (v) to (-0.16,0.34409);         
    \draw (v) to (-0.16, -0.34409);

    % Vertices u_ks
    \node[node, label=above:{$u_i$}] (ui) at (1,0.25) {};
    \node[node, label=below:{$u_j$}] (uj) at (1,-0.25) {};
    \node[node] (u1) at (0.7,0.5) {};
    \node[node] (u2) at (0.7,-0.5) {};

    % edges between u_ks
    \draw (u1) -- (u2) -- (uj) -- (ui) -- (u2);
    \draw (ui) -- (u1) -- (uj);

    % Edges between v and u_ks
    \draw[bend left=10] (v) to (u1);
    \draw (v) -- (u1);
    \draw[bend right=10] (v) to (u1);
    \draw[bend left=10] (v) to (u2);
    \draw (v) -- (u2);
    \draw[bend right=10] (v) to (u2);
    \draw[bend left=10] (v) to (ui);
    \draw (v) -- (ui);
    \draw[bend right=10] (v) to (ui);
    \draw[bend left=10] (v) to (uj);
    \draw (v) -- (uj);
    \draw[bend right=10] (v) to (uj);

    % Branches of the unnamed vertices
    \draw (0.59,0.79) circle [radius=0.13];
    \draw (0.59,-0.79) circle [radius=0.13];
    \draw[bend right=8] (u1) to (0.718, 0.767284);
    \draw[bend left=8] (u1) to (0.718, 0.767284);
    \draw (u1) to (0.496,0.7002);
    \draw[bend right=8] (u2) to (0.496,-0.7002);
    \draw[bend left=8] (u2) to (0.496,-0.7002);
    \draw (u2) to (0.718, -0.767284);
    
    % B_i
    \draw (1.7,0.4) ellipse [x radius=0.3, y radius=0.23];
    \node at (2.06,0.58) {$B_i$};
    \node[node, label=right:{$b_i$}] (bi) at (1.64,0.4) {};
    \draw[bend right=8] (ui) to (1.426, 0.493659);
    \draw[bend left=8] (ui) to (1.426, 0.493659);
    \draw (ui) to (1.408,0.347239);
    \draw (ui) to (1.495, 0.23208);
    
    % B_j
    \draw (1.7,-0.4) ellipse [x radius=0.3, y radius=0.23];
    \node at (2.06,-0.6) {$B_j$};
    \node[node, label=right:{$b_j$}] (bj) at (1.64,-0.4) {};
    \draw (uj) to (1.426, -0.493659);
    \draw (uj) to (1.408,-0.347239);
    \draw (uj) to (1.495, -0.23208);
    
\end{tikzpicture}
\caption{Sketch of $G$ with two non-empty branches.}
\label{fig:7}
\end{figure}

\begin{claim}\label{claim:oneNonEmptyBranch}
    There is up to one non-empty branch.
\end{claim}
\begin{proof}
    Assume for a contradiction there are two different non-empty branches $B_i$ and $B_j$ attached to $u_i$ and $u_j$ respectively. So there exist vertices $b_i\in B_i$ and $b_j \in B_j$, as sketched in \cref{fig:7}. We will show that the divisor $D_2-b_i-b_j$ has rank at least $1$, which is a contradiction since its degree is smaller than the gonality. 

By \cref{claim:kChipsOnv}, $D_2-b_i-b_j-w$ is equivalent to an effective divisor for any $w\in V(G'')$. It suffices to show that $D_2-b_i-b_j$ is equivalent to effective when taking away a chip from either $b_i$ or $b_j$ as there is already a chip everywhere else.

Consider first the divisor $D_2 - 2b_i$. Performing Dhar's burning algorithm at $b_i$ results in an effective divisor since the divisor has rank at least $0$. If the algorithm results in any changes to the divisor on $B_j$, there is an iteration in which the fire burns some, but not all of the vertices of  $B_j\cup\{u_j\}$. By \cref{claim:branches}, this can never happen and the configuration of chips on $B_j$ does not change during the  algorithm. Since $D_2$ had a chip at $b_j$, the resulting divisor has a chip at $b_j$ as well. 

It follows that the divisor $D_2 - 2b_i - b_j$ is equivalent to an effective divisor. 
By a similar argument, $D_2 - b_i - 2b_j$ is equivalent to an effective divisor as well, so we are done. 
\end{proof}

\begin{figure}
\centering
\begin{tikzpicture}[scale=2, node/.style={circle, draw, fill=black, inner sep=1.2pt}]
    % Vertices for first section
    \node[node, label=below left:{$v$}] (v) at (0,0) {};
    \node[node, label=below:{$u_1=u$}] (u1) at (2.2,0) {};
    \node[node, label=right:{$u_2$}] (u2) at (2.01262,0.84524) {};
    \node[node, label=above:{$u_3$}] (u3) at (1.48558,1.53209) {};
    \node[node, label=left:{$u_4$}] (u4) at (0.71764,1.93185) {};
    
    % Edges for first section
    \draw (u1)--(u2)--(u3)--(u4)--(u1)--(u3)--(u1);
    \draw (u2)--(u4);
    \draw[bend left=6] (v) to (u1);
    \draw[bend right=6] (v) to (u1);
    \draw (v) to (u1);
    \draw[bend left=6] (v) to (u2);
    \draw[bend right=6] (v) to (u2);
    \draw (v) to (u2);
    \draw[bend left=6] (v) to (u3);
    \draw[bend right=6] (v) to (u3);
    \draw (v) to (u3);
    \draw[bend left=6] (v) to (u4);
    \draw[bend right=6] (v) to (u4);
    \draw (v) to (u4);

    % Drawing G''
    \draw (-1.6,0.7) circle [radius=0.7];
    \node at (-1.6,0.7) {$G''$};
    \draw (v) to (-1.6, 0);
    \draw (v) to (-1.06, 1.14542);
    \draw (v) to (-0.92, 0.53387);
    \draw[bend right=7] (v) to (-1.086, 0.22481);
    \draw[bend left=7] (v) to (-1.086, 0.22481);

    % Drawing B
    \draw (3.8,0.7) circle [radius=0.7];
    \node at (3.8,0.7) {$B$};
    \draw (u1) to (3.8, 0);
    \draw (u1) to (3.26, 1.14542);
    \draw[bend right=7] (u1) to (3.18,0.37504);
    \draw[bend left=7] (u1) to (3.18,0.37504);
    
\end{tikzpicture}
\caption{Structure of $G$. Both $G''$ and $B$ are (possibly disconnected or empty) subgraphs of $G$. There are $k + 1$ edges between $v$ and each $u_i$ and one edge between each distinct pair of $u_i$.}
\label{fig:8}
\end{figure}
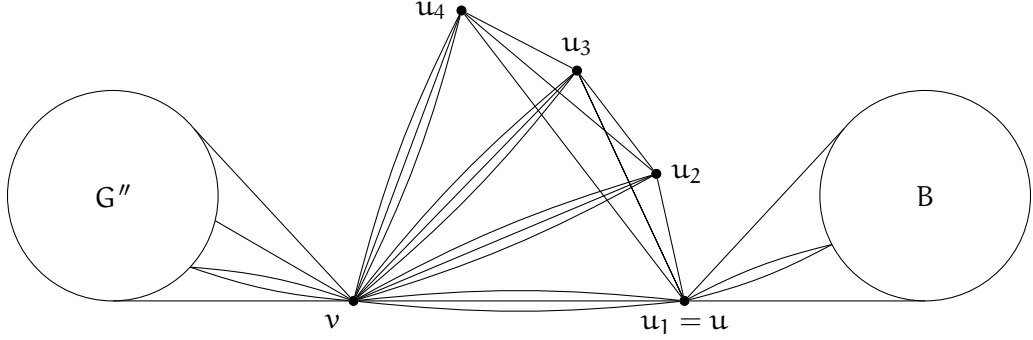

Let $B:=G'\setminus N_v(G')$.
Note that the subgraph $B$ is possibly disconnected, or empty.
By \cref{claim:WIsEmpty,claim:oneNonEmptyBranch}, if $V(B)$ is non-empty, it is equal to a single branch $B_i$ of some $u_i$. We relabel that vertex to $u$.
The existence of $D_2$ severely restricts the structure of the graph $G$, as seen in \cref{fig:8}.

\begin{claim}\label{claim:divisorD_2Prime}
    Using the notation from \cref{fig:8}, the $u$-reduced divisor equivalent to $D_2$ is
    \[
    D_2' :=(k+n+1)\cdot u + \sum_{w\in V(B)} w.
    \]
\end{claim}
\begin{proof}
    We apply Dhar's burning algorithm to $D_2$, starting a fire at $u$. By \cref{claim:branches}, all of $B$ burns down, while no other vertices do. Now set-firing the unburnt vertices, namely 
    \[
    V(G'')\cup\{v\}\cup N_{G'}(v)\setminus\{u\},
    \] 
    gives 
    \[
    (k+n+1)\cdot u + \sum_{w\in V(B)} w,
    \]
    which is $u$-reduced by \cref{claim:branches}.
\end{proof}

This completes the proof of \cref{lemma:structureMainLemma}. Moreover, when $B$ is non-empty we can iterate the argument by taking $u$ as the new vertex $v$, the divisor $D_2'$ as $D_2$, $B$ as the subgraph $G'$ and $G''\cup\big(\{v\}\cup N_v(G')\setminus\{u\}\big)$ as $G''$, but more on that in \cref{proposition:firstDirectionOfMain}.

\subsection{Quasi-banana graphs}\label{sec:quasiBananaGraphs}

In this subsection, we construct a large family of gonality-tight graphs. We show that this family includes all gonality-tight graphs $G$ whose second gonality is realised by a divisor of the form
\begin{equation*}
    D=v+\sum_{w\in V(G)} w.
\end{equation*}
We refer to these graphs as \emph{quasi-banana graphs}.

\begin{definition}\label{def:quasiBananaGraph} 
    Let $m$ be a non-negative integer and let $N=(q_1,q_2,\ldots, q_m)$  be an $m$-tuple of positive integers. 
    A \emph{quasi-banana graph}  $QB_m(N)$
is defined as follows. The vertex set is 
    \[
    V\left(QB_m(N)\right)=\bigcup_{i=0}^m C_i,
    \] 
    where,
    \[
    C_0=\{u_{0,1}\}\qquad\text{and} \qquad C_i=\{u_{i,j} \mid  1\leq j\leq q_i\}.
    \]
    The edge set is specified by the following.
    \begin{enumerate}
    \item For every fixed $1\leq i\leq m$, and distinct $1\leq j,k\leq q_i$, there is exactly one edge between $u_{i,j}$ and $u_{i,k}$.
    \item For every fixed $0\leq i<m$, and $1\leq j\leq q_{i+1}$, there are exactly $\sum_{k=1}^iq_k+2$ edges between $u_{i,1}$ and $u_{i+1,j}$.
    \item There are no other edges.
    \end{enumerate}
\end{definition}

\begin{example}
    The quasi-banana graph $QB_0$ consists of a single vertex and no edges.
\end{example}

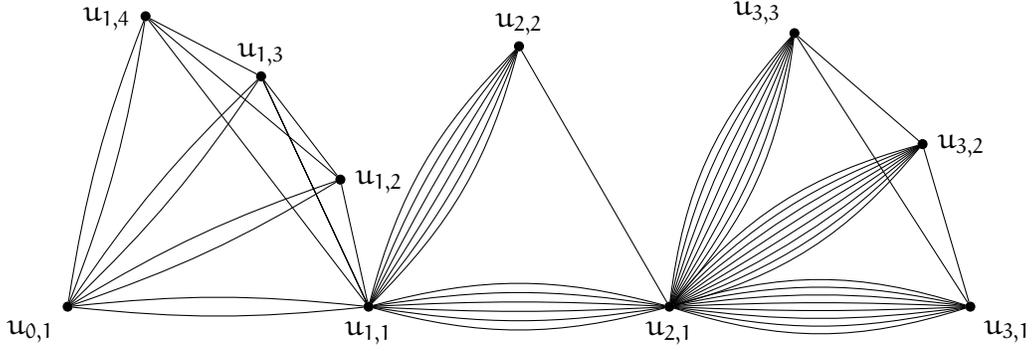
\begin{figure}
\centering
\begin{tikzpicture}[scale=2, node/.style={circle, draw, fill=black, inner sep=1.2pt}]
    % Vertices for first section
    \node[node, label=below left:{$u_{0,1}$}] (v) at (0,0) {};
    \node[node, label=below:{$u_{1,1}$}] (u11) at (2,0) {};
    \node[node, label=right:{$u_{1,2}$}] (u12) at (1.81262,0.84524) {};
    \node[node, label=above:{$u_{1,3}$}] (u13) at (1.28558,1.53209) {};
    \node[node, label=left:{$u_{1,4}$}] (u14) at (0.51764,1.93185) {};
    
    % Edges for first section
    \draw (u11)--(u12)--(u13)--(u14)--(u11)--(u13)--(u11);
    \draw (u12)--(u14);
    \draw[bend left=6] (v) to (u11);
    \draw[bend right=6] (v) to (u11);
    \draw[bend left=6] (v) to (u12);
    \draw[bend right=6] (v) to (u12);
    \draw[bend left=6] (v) to (u13);
    \draw[bend right=6] (v) to (u13);
    \draw[bend left=6] (v) to (u14);
    \draw[bend right=6] (v) to (u14);

    % Vertices for second section
    \node[node, label=below:{$u_{2,1}$}] (u21) at (4,0) {};
    \node[node, label=above:{$u_{2,2}$}] (u22) at (3,1.73205) {};

    % Edges for second section
    \draw (u21)--(u22);
    \draw[bend left=3]  (u11) to (u21);
    \draw[bend right=3] (u11) to (u21);
    \draw[bend left=9] (u11) to (u21);
    \draw[bend right=9](u11) to (u21);
    \draw[bend left=15] (u11) to (u21);
    \draw[bend right=15](u11) to (u21);
    \draw[bend left=3]  (u11) to (u22);
    \draw[bend right=3] (u11) to (u22);
    \draw[bend left=9] (u11) to (u22);
    \draw[bend right=9](u11) to (u22);
    \draw[bend left=15] (u11) to (u22);
    \draw[bend right=15](u11) to (u22);

    % Vertices for third section
    \node[node, label=below right:{$u_{3,1}$}] (u31) at (6,0) {};
    \node[node, label=right:{$u_{3,2}$}] (u32) at (5.68251,1.08128) {};
    \node[node, label=above left:{$u_{3,3}$}] (u33) at (4.83083,1.81926) {};

    % Edges for third section
    \draw (u31)--(u32)--(u33)--(u31);
    \draw[bend left=2.5]  (u21) to (u31);
    \draw[bend right=2.5] (u21) to (u31);
    \draw[bend left=7.5]  (u21) to (u31);
    \draw[bend right=7.5] (u21) to (u31);
    \draw[bend left=12.5]  (u21) to (u31);
    \draw[bend right=12.5] (u21) to (u31);
    \draw[bend left=17.5]  (u21) to (u31);
    \draw[bend right=17.5] (u21) to (u31);
    \draw[bend left=2.5]  (u21) to (u32);
    \draw[bend right=2.5] (u21) to (u32);
    \draw[bend left=7.5]  (u21) to (u32);
    \draw[bend right=7.5] (u21) to (u32);
    \draw[bend left=12.5]  (u21) to (u32);
    \draw[bend right=12.5] (u21) to (u32);
    \draw[bend left=17.5]  (u21) to (u32);
    \draw[bend right=17.5] (u21) to (u32);
    \draw[bend left=2.5]  (u21) to (u33);
    \draw[bend right=2.5] (u21) to (u33);
    \draw[bend left=7.5]  (u21) to (u33);
    \draw[bend right=7.5] (u21) to (u33);
    \draw[bend left=12.5]  (u21) to (u33);
    \draw[bend right=12.5] (u21) to (u33);
    \draw[bend left=17.5]  (u21) to (u33);
    \draw[bend right=17.5] (u21) to (u33);
\end{tikzpicture}
\caption{The quasi-banana graph $QB_3(4,2,3)$.}
\label{fig:QB_3(423)}
\end{figure}

\begin{example}
    The quasi-banana graph $QB_3(4,2,3)$ is illustrated in \cref{fig:QB_3(423)}.
\end{example}

\begin{example}

The \emph{generalised banana graphs} $B_{m+1,m+1}^{*}$, introduced in \cite{ADMYY_GonalitySequencesOfGraphs} (although initially studied in \cite{CastryckCools_NewtonAndGonalities} by a different name)  
and further studied in \cite{FJKO_semigroupOfGonality},
coincide with the  quasi-banana graphs $QB_{m}(\mathbf{1}_{m})$, where $\mathbf{1}_{m}$ is the $m$-tuple with all entries equal to one.

\end{example}

\begin{example}

    Let $QB_m(q_1,q_2,\ldots, q_m)$ be a quasi-banana graph with $m\geq 1$. Then 
    \[
    QB_m(q_1,q_2,\ldots, q_m)\setminus C_m \cong QB_{m-1}(q_1,q_2,\ldots, q_{m-1}).
    \]
    If $q_m\geq2$ then
    \[
    QB_m(q_1,q_2,\ldots, q_m)\setminus \{u_{m,q_m}\} \cong QB_{m}(q_1,q_2,\ldots, q_m-1).
    \]
\end{example}

Before describing the connection to gonality-tight graphs, we record the following claim.

\begin{claim}\label{claim:everyFireBurns}
Let $H$ be a graph whose first gonality equals $|V(H)|$ and fix a vertex $v$ of $H$. Let $E$ be an effective divisor with at most a single chip on every vertex other than $v$. Then  applying the Dhar's burning algorithm from $v$ results in burning the entire graph. Namely, $E$ is $v$-reduced.  
\end{claim}

\begin{proof}
   The claim doesn't depend on the number of chips at $v$, so we may assume that $E(v)=1$. We may also assume that $E$ has a chip at every vertex. Suppose that the fire doesn't burn the entire graph. Then there is a collection  of vertices $\emptyset\neq U\neq V(H)$ such that $E+\ddiv(U)$  in an effective divisor. But then there is a vertex $w\notin U$ such that $\left(E+\ddiv(U)\right)(w) \geq 1 $. 
   It follows that the divisor $E-w$ has rank $1$, which is a contradiction, since its degree is smaller than $|V(H)|$. 
\end{proof}

The following proposition is the main result of the section.

\begin{proposition}\label{proposition:firstDirectionOfMain}
    Let $G$ be a gonality-tight graph with a vertex $v$ such that the divisor
    \[
    D=v+\sum_{w\in V(G)} w
    \]
    realises the second gonality. Then $G$ is a quasi-banana graph and $v$ is the vertex $u_{0,1}$ in the notation of \cref{def:quasiBananaGraph}.
\end{proposition}

\begin{proof} 
Since $G$ is gonality-tight, and
\[
\gon_2(G)=\deg(D)=|V(G)|+1,
\]
it follows that $\gon_{1}(G)=|V(G)|$. By \cref{claim:everyFireBurns}, the divisor $D$ is $v$-reduced. Relabel the vertex $v$ as $u_{0,1}$, and set $C_0=\{u_{0,1}\}$. Then
\[
D=2u_{0,1}+\sum_{w\in V(G)\setminus C_0}w.
\]
We now apply \cref{lemma:structureMainLemma} to $G$ with $k=1$, taking $u_{0,1}$ as the vertex $v$, $G'=G\setminus \{u_{0,1}\}$ and $G''=\emptyset$. 

According to the lemma, if $V(G)\setminus \{u_{0,1}\}$ is non-empty, there are exactly $k+1 = 2$  edges between $u_{0,1}$ and each of its neighbours in $G'$. Denote the set of neighbours of $u_{0,1}$ by 
\[
C_1=N_{G'}(u_{0,1})=\{u_{1,1},u_{1,2},\ldots, u_{1,q_1}\},
\]
taking the vertex $u_{1,1}$ as the unique vertex $u$ from Part (3) of the lemma if $V(G)\setminus(C_0\cup C_1)$ is non-empty. The vertices in $C_1$ form a complete subgraph, and each one is connected to $u_{0,1}$ by exactly two edges. 
Part (4) then implies that the $u_{1,1}$-reduced divisor equivalent to $D$ is
\[
D'=(q_1+2)\cdot u_{1,1}+\sum_{w\in V(G)\setminus(C_0\cup C_1)}w. %\sim D.
\]

If $V(G)\setminus(C_0\cup C_1)$ is non-empty, we apply \cref{lemma:structureMainLemma} again, taking $u_{1,1}$ as the new vertex $v$, and letting $G'=G\setminus(C_0\cup C_1)$.
The vertices in the next set of neighbours
\[
C_2=N_{G'}(u_{1,1})=\{u_{2,1},u_{2,2},\ldots, u_{2,q_2}\}
\]
are connected with $q_1+2$ edges to $u_{1,1}$ each and have precisely one edge between each other. If $V(G)\setminus (C_0\cup C_1\cup C_2)$ is non-empty, Part (3) of the lemma again provides a unique vertex $u_{2,1}\in C_2$ separating the remaining vertices.

Continuing inductively, after $i$ applications of the lemma, we have constructed disjoint vertex sets 
\[
C_0,C_1,\ldots, C_{i},
\]
where for each $j\geq1$ we write $C_j=\{u_{j,1},u_{j,2},\ldots, u_{j,q_j}\}$. 
These satisfy the following.
\begin{itemize}
    \item For each $j\geq 1$ there is exactly one edge between each distinct pair of vertices in $C_j$.
    \item For $j\geq1$, each vertex in $C_j$ is connected to $u_{j-1,1}$ with $\sum_{k=1}^{j-1}q_k+2$ edges, where $q_k=|C_k|$.
    \item If $ V(G)\setminus \bigcup_{k=1}^iC_k$ is non-empty, then the vertex $u_{i,1}$ is a cut vertex separating the remaining vertices. The $u_{i,1}$-reduced divisor realising the second gonality is
    \[
    \left(\sum_{k=1}^iq_k+2\right)\cdot u_{i,1}+\sum_{w\in V(G)\setminus \bigcup_{k=1}^iC_k}w\sim D.
    \]
    This allows us to apply the lemma again.
\end{itemize}
As $G$ is finite, this process eventually terminates after $m$ iterations, for some $m\geq0$. 

The resulting graph is completely determined by the $q_i=|C_i|$, and satisfies the conditions of \cref{def:quasiBananaGraph}. Hence $G$ is a quasi-banana graph $QB_m(q_1,q_2,\ldots, q_m)$, with $u_{0,1}=v$.
\end{proof}

We will show in 
\cref{sec:proofOfMainTheorem} that every quasi-banana graph $G$ is gonality-tight and prove \cref{mainTheorem:alternativeDefinitonForQBGraphs}. At this point, we have only established that the family of gonality-tight graphs with a divisor of the form $D$ realising the second gonality is a subfamily of quasi-banana graphs.

\section{Two technical lemmas} \label{sec:TwoTechnicalLemmas}%\simun{I think this title sounds good.}}

In this section, we prove two technical results, namely \cref{lemma:upperBoundMainLemma,lemma:lowerBoundMainLemma} that will be used inductively in \cref{sec:proofOfMainTheorem} to prove the upper and lower bound respectively of the gonality sequence.

\subsection{Bounding gonalities from above}
Recall that, when $D$ is a divisor on a subgraph $H$ of $G$, we denote by $r_H(D)$ and $r_G(D)$ the rank of $D$ as a divisor on $H$ and $G$ respectively. Similarly, we denote the genus of a graph $H$ by $g_H$.
This subsection is devoted to  the following lemma.

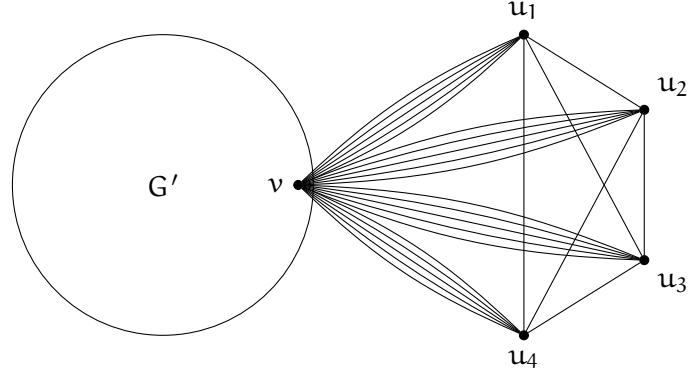
\begin{figure}
\centering
\begin{tikzpicture}[scale=2, node/.style={circle, draw, fill=black, inner sep=1.2pt}]
    % Vertex v
    \node[node, label=left:{$v$}] (v) at (0.3,0) {};

    % Graph G'
    \draw (-0.6,0) circle [radius=1];
    \node at (-0.6,0) {$G'$};

    % Complete graph K_n
    %\node at (2.3,1.1) {$K_n$};
    \node[node, label=above:{$u_1$}] (u1) at (1.8,1) {};
    \node[node, label=above right:{$u_2$}] (u2) at (2.6,0.5) {};
    \node[node, label=below right:{$u_3$}] (u3) at (2.6,-0.5) {};
    \node[node, label=below:{$u_4$}] (u4) at (1.8,-1) {};

    % Edges of K_n
    \draw (u1) -- (u2);
    \draw (u1) -- (u3);
    \draw (u1) -- (u4);
    \draw (u2) -- (u3);
    \draw (u2) -- (u4);
    \draw (u3) -- (u4);

    % Edges between v and K_n
    \draw[bend right=0](v) to (u1);
    \draw[bend left=5]  (v) to (u1);
    \draw[bend right=5] (v) to (u1);
    \draw[bend left=10] (v) to (u1);
    \draw[bend right=10](v) to (u1);

    \draw[bend right=0](v) to (u2);
    \draw[bend left=5]  (v) to (u2);
    \draw[bend right=5] (v) to (u2);
    \draw[bend left=10] (v) to (u2);
    \draw[bend right=10](v) to (u2);
    
    \draw[bend right=0](v) to (u3);
    \draw[bend left=5]  (v) to (u3);
    \draw[bend right=5] (v) to (u3);
    \draw[bend left=10] (v) to (u3);
    \draw[bend right=10](v) to (u3);

    \draw[bend right=0](v) to (u4);
    \draw[bend left=5]  (v) to (u4);
    \draw[bend right=5] (v) to (u4);
    \draw[bend left=10] (v) to (u4);
    \draw[bend right=10](v) to (u4);

\end{tikzpicture}

\caption{Graph $G$. Here $k=4$ and $n=4$.}
\label{fig:G}
\end{figure}

\begin{lemma}\label{lemma:upperBoundMainLemma}
    Let $k$ be a non-negative 
    integer and suppose that a graph $G'$ of genus $g_{G'}=\binom{k}{2}$ has a vertex $v$ such that, for every integer $1\leq l\leq k-2$,
    \begin{equation*}
        r_{G'}\big(l (k+1)\cdot v\big) \geq  \frac{l(l+3)}{2}.
    \end{equation*}
    Let $G$ be the graph obtained by attaching a complete graph $K_n$ to $G'$ along 
$k+1$ edges between $v$ and each vertex in $K_n$ (see \cref{fig:G} for an illustration). Then the following hold. 
    \begin{enumerate}
        \item The genus of $G$ is $g_G=\binom{k+n}{2}$.
        \item For every $u\in V(K_n)$, the divisor $(k+n-2)(k+n+1)\cdot u$ is equivalent to the canonical divisor $K_G$.
        \item For every $-2\leq l\leq k+n$  and every $u\in V(K_n),$
        \begin{equation*}
              r_{G}\big(l(k+n+1)\cdot u\big) \geq  \frac{l(l+3)}{2}.
        \end{equation*}
    \end{enumerate}
\end{lemma}

\medskip

\subsubsection*{Proof of \cref{lemma:upperBoundMainLemma}(1)} 
The genus of $G$ can be computed explicitly from the definition: 
\begin{align*}
    g_G&=|E(G)|-|V(G)|+1 \\
    &=\left(|E(G')|+n(k+1)+\frac{n(n-1)}{2}\right)-\left(|V(G')|+n\right)+1 \\
    &=\left(|E(G')|-|V(G')|+1\right)+n(k+1)+\frac{n(n-1)}{2}-n \\
    &=g_{G'}+\frac{2nk+n^2-n}{2}
    =\frac{k^2-k+2nk+n^2-n}{2} \\
    &=\frac{(k+n)^2-(k+n)}{2}
    =\binom{k+n}{2}.
\end{align*}

\subsubsection*{Proof of \cref{lemma:upperBoundMainLemma}(2)}

Now, consider the divisor
\[
D_2 := (k+1)\cdot v + \sum_{w\in V(K_n)} w.
\]

\begin{claim}\label{claim:divisorD_2}
$D_2$ satisfies the following:
\begin{enumerate}
    \item For every $u\in V(K_n)$, $D_2\sim (k+n+1)\cdot u$.
    \item $(k+1)\cdot D_2\sim (k+n+1)(k+1)\cdot v$.
\end{enumerate}
\end{claim}
\begin{proof}
Both statements are consequences of set-firing.
Set-firing $V(G)\setminus\{u\}$ with respect to $D_2$ gives $(k+n+1)\cdot u$, and set-firing $V(K_n)$ with respect to $(k+1)\cdot D_2$ produces $(k+n+1)(k+1)\cdot v$.
\end{proof}

Returning to prove \cref{lemma:upperBoundMainLemma}(2), notice that
 \[
 r_{G'}\big((k-2)(k+1)\cdot v\big)\geq\frac{(k-2)(k+1)}{2}= \binom{k}{2}-1=g_G-1
 \]
and
\[
\deg\big((k-2)(k+1)\cdot v\big) = (k-2)(k+1) = k(k-1)-2= 2g_{G'}-2.
\]
By \cref{corollary:corollariesOfRiemann-Roch}, we have  $(k-2)(k+1)\cdot v\sim K_{G'}$.
% \simun{I added the following sentence.} \yoav{Sure. Is this notation only used here or elsewhere as well?} \simun{This is the only time we use this notation. }
In the following argument, whenever $u$ is a vertex in a subgraph $H'$ of $H$, we use $\operatorname{val}_{H'}(u)$ to denote the valency of $u$ considered as a vertex in $H'$.
The canonical divisor of $G$ can now be computed:
\begin{align*}
    K_G&=\sum_{w\in V(G)}\big(\operatorname{val}_{G}(w)-2\big)\cdot w = \sum_{w\in V(G')}\big(\operatorname{val}_{G}(w)-2\big)\cdot w + \sum_{w\in V(K_n)}\big(\operatorname{val}_{G}(w)-2\big)\cdot w \\
    &=\sum_{w\in V(G')\setminus\{v\}}\big(\operatorname{val}_{G'}(w)-2\big)\cdot w + \big(\operatorname{val}_{G'}(v)+n(k+1)-2\big)\cdot v + \sum_{w\in V(K_n)}\big(\operatorname{val}_{G}(w)-2\big)\cdot w \\
    &=\sum_{w\in V(G')}\big(\operatorname{val}_{G'}(w)-2\big)\cdot w +n(k+1)\cdot v +\sum_{w\in V(K_n)}\big((k+n)-2\big)\cdot w \\
    &=K_{G'}+n(k+1)\cdot v +\sum_{w\in V(K_n)}(k+n-2)\cdot w \\
    &\sim(k+1)(k-2)\cdot v + n(k+1)\cdot v +\sum_{w\in V(K_n)}(k+n-2)\cdot w \\
    %&=(k+n-2)(k+1)\cdot v + \sum_{w\in V(K_n)}(k+n-2)\cdot w \\
    &=(k+n-2)\cdot D_2
\end{align*}
By \cref{claim:divisorD_2}, for any $u\in V(K_n)$, 
\[
K_G\sim (k+n-2)\cdot D_2\sim (k+n-2)(k+n+1)\cdot u.
\]
This concludes the proof of Part (2).

\subsubsection*{Proof of \cref{lemma:upperBoundMainLemma}(3)}
This is the trickiest and most technical part of the proof of the lemma. 
We  first show that, for every $u\in V(K_n)$ and every  integer $-2\leq l\leq \frac{k+n-2}{2}$, 
\begin{equation*}
    r_{G}\big(l(k+n+1)\cdot u\big)=r_{G}(l\cdot D_2) \geq  \frac{l(l+3)}{2}.
\end{equation*}
The first equality follows from \cref{claim:divisorD_2}, and we prove the inequality by strong induction on $l$. For the base cases $l = -2, -1,$ or $0$, the statement holds since 
\[
r_G(-2 \cdot D_2) = r_G(-1 \cdot D_2) = -1 \quad \text{and} \quad r_G(0 \cdot D_2) = 0.
\]
Now fix some $1\leq l\leq \frac{k+n-2}{2}$, and assume that
\[
r_G(r\cdot D_2)\geq \frac{r(r+3)}{2} \qquad\text{for all } -2 \leq r \leq l-1.
\]
Let $E$ be any effective divisor on $G$ with $\deg(E)=\frac{l(l+3)}{2}$. We want to show that $l\cdot D_2 - E$ is equivalent to an effective divisor. First, suppose $1\leq l\leq k$. We have the following two cases.

\casetitle[$l\cdot D_2 - E$ is not effective on $K_n$]{A}

\medskip
\noindent
Here, there must exist $w\in V(K_n)$, such that $(l+1)\cdot w\leq E$. Therefore, it suffices to show that
\[
r_G\big(l\cdot D_2-(l+1)\cdot w\big)\geq \frac{l(l+3)}{2}-(l+1)=\frac{(l-1)(l+2)}{2}.
\]
By \cref{claim:divisorD_2},
\[
l\cdot D_2-(l+1)\cdot w\sim (l-1)\cdot D_2+\big(k+n+1-(l+1)\big)\cdot w\geq (l-1)\cdot D_2.
\]
Hence, by the inductive hypothesis,
\[
r_G\big(l\cdot D_2-(l+1)\cdot w\big)\geq r_G\big((l-1)\cdot D_2\big)\geq \frac{(l-1)(l+2)}{2}.
\]
This concludes Case (A).

\casetitle[$l\cdot D_2 - E$ is effective on $K_n$]{B}

\medskip
\noindent
If $l\cdot D_2 - E$ is effective on $K_n$, then by \cref{lemma:restrictedRank}, it suffices to show that
\[
r_{G'}(l\cdot D_2|_{G'})=r_{G'}\big(l(k+1)\cdot v\big)\geq \frac{l(l+3)}{2}.
\]
This follows from the following claim.

\begin{claim}\label{claim:rankOnG'}
    \[
    r_{G'}\big(r(k+1)\cdot v\big)\geq
    \begin{cases}
        \frac{r(r+3)}{2} & \text{if } 0\leq r\leq k,\\
        \frac{r(r+3)}{2}-\frac{(r-k)(r-k+1)}{2}  & \text{if } k+1\leq r.
    \end{cases}
    \]
\end{claim}
\begin{proof}
    If $r=0$, then clearly $r_{G'}\big(0(k+1)\cdot v\big)=0$.
    For $1\leq r\leq k-2$, the inequality follows from the assumption. When $r\geq k-1$, we have
    \[
    \deg\big(r(k+1)\cdot v\big)=r(k+1)\geq (k-1)(k+1)=k(k-1)+k-1\geq 2g_{G'}-1.
    \]
    Thus, by \cref{corollary:corollariesOfRiemann-Roch},
    \begin{align*}
        r_{G'}\big(r(k+1)\cdot v\big)&=\deg\big(r(k+1)\cdot v\big)-g_{G'}=r(k+1)-\frac{k(k-1)}{2} \\
        %&=\frac{2rk+2r-k^2+k}{2}=\frac{(r^2+3r)-(r^2-rk+r-kr+k^2-k)}{2} \\
        %&=\frac{r(r+3)}{2}-\frac{r(r-k+1)-k(r-k+1)}{2} \\
        &=\frac{r(r+3)}{2}-\frac{(r-k)(r-k+1)}{2}.
    \end{align*}
    Note that the second term is $0$ when $r=k-1$ or $r=k$. This proves the claim.
\end{proof}

Now, assume that $k+1\leq l\leq \frac{k+n-2}{2}$. This time, there are 3 cases to consider:

\casetitle[$l\cdot D_2 - E$ is not effective on $K_n$]{C}

\medskip
\noindent
The argument is the same as before. There exists $w\in V(K_n)$ such that $(l+1)\cdot w\leq E$, and therefore
\[
r_G\big(l\cdot D_2-(l+1)\cdot w\big)\geq r_G\big((l-1)\cdot D_2\big)\geq \frac{(l-1)(l+2)}{2}= \frac{l(l+3)}{2}-(l+1),
\]
as in the previous Case (A).

\casetitle[$l\cdot D_2 - E$ is effective on $K_n$, and $\deg(E|_{K_n}) \geq \frac{(l-k)(l-k+1)}{2}$]{D}

\medskip
\noindent
As $l\cdot D_2-E|_{K_n}$ is effective, by \cref{lemma:restrictedRank}, it is enough to show that $r_{G'}\big((l\cdot D_2-E|_{K_n})|_{G'}\big)\geq \deg(E|_{G'})$. From \cref{claim:rankOnG'}, we have
\begin{align*}
    r_{G'}\big((l\cdot D_2-E|_{K_n})|_{G'}\big)&=r_{G'}\big(l(k+1)\cdot v\big)\geq\frac{l(l+3)}{2}-\frac{(l-k)(l-k+1)}{2} \\
    &\geq \deg(E)-\deg(E|_{K_n})=\deg(E|_{G'}).
\end{align*}
The desired inequality holds, completing Case (D).

\casetitle[$l\cdot D_2 - E$ is effective on $K_n$, and $\deg(E|_{K_n}) < \frac{(l-k)(l-k+1)}{2}$]{E}

\medskip
\noindent
Let $0\leq r_0\leq l-k-1$ be the unique integer such that
\begin{equation}\label{eq:r_0Inequality}
    \frac{r_0(r_0+1)}{2}\leq \deg(E|_{K_n})\leq \frac{(r_0+1)(r_0+2)}{2}-1=\frac{r_0(r_0+3)}{2}.
\end{equation}
Next, set 
\[
q=\left\lfloor \frac{l-r_0}{k+1} \right\rfloor
\]
and let $0\leq s\leq k$ be the remainder of the division of $l-r_0$ by $k+1$. So,
\[
l-r_0=q(k+1)+s.
\]
Note that $q\geq1$.
By \cref{claim:divisorD_2}, 
\begin{equation}\label{eq:rD_2Equivalence}
    \begin{aligned}
        l\cdot D_2%&= r_0\cdot D_2 + (l-r_0)\cdot D_2\\ 
        &= r_0\cdot D_2+q(k+1)\cdot  D_2 + s\cdot D_2 \\
        &\sim (r_0+s)\cdot D_2+q(k+n+1)(k+1)\cdot v
    \end{aligned}
\end{equation}
Since $r_0+s = l-q(k+1)<l$, the inductive hypothesis gives,
\begin{equation*}
    r_{G}\big((r_0+s)\cdot D_2\big)\geq \frac{(r_0+s)(r_0+s+3)}{2}\geq \frac{r_0(r_0+3)}{2}\geq \deg(E|_{K_n})
\end{equation*}
Therefore, by \cref{eq:r_0Inequality},
\begin{equation}\label{eq:rankOfDivisorIsA}
\begin{aligned}
    r_G\big((r_0+s)\cdot D_2-E|_{K_n}\big) &\geq \frac{(r_0+s)(r_0+s+3)}{2} - \deg(E|_{K_n}) \\
    &\geq \frac{(r_0+s)(r_0+s+3)}{2}-\frac{r_0(r_0+3)}{2} =: \gamma.
\end{aligned}
\end{equation}
Write $E|_{G'}$ as a sum of two effective divisors $E_1$ and $E_2$, where $\deg(E_1)=\max\big(\gamma,\deg(E|_{G'})\big)$. 

\begin{claim}\label{claim:1stEffectiveDivisor}
    The divisor $(r_0+s)\cdot D_2-E|_{K_n}-E_1$ is equivalent to an effective divisor on $G$.
\end{claim}
\begin{proof}
    This follows immediately from \cref{eq:rankOfDivisorIsA} and the fact that $\deg(E_1)\leq \gamma$.
\end{proof}

\begin{claim}\label{claim:2ndEffectiveDivisor}
    The divisor $q(k+n+1)(k+1)\cdot v-E_2$ is equivalent to an effective divisor on $G$.
\end{claim}
\begin{proof}
If $\deg(E_1)=\deg(E|_{G'})$, then $E_2$ is the zero divisor, and the claim holds. So, we can now assume that $\deg(E_1)=\gamma$. From \cref{eq:r_0Inequality},
\begin{equation}
\begin{aligned}\label{eq:degreeInequality}
    \deg(E_2)&=\deg(E|_{G'})-\deg(E_1)=\deg(E)-\deg(E_1)-\deg(E|_{K_n}) \\
    &\leq \frac{l(l+3)}{2}-\gamma-\frac{r_0(r_0+1)}{2}.
\end{aligned}
\end{equation}
As $E_2$ is supported on $G'$, by \cref{lemma:restrictedRank}, it suffices to show that:
\begin{equation*}
    r_{G'}\big(q(k+n+1)(k+1)\cdot v\big)\geq \deg(E_2).
\end{equation*}
Since $q\geq1$, by \cref{claim:rankOnG'},
\begin{equation}\label{eq:rankInequality}
r_{G'}\big(q(k+n+1)(k+1)\cdot v\big)\geq \frac{\alpha(\alpha+3)}{2} \\ -\frac{(\alpha-k)(\alpha-k+1)}{2},
\end{equation}
where $\alpha:=q(k+n+1)$.

Combining \cref{eq:degreeInequality,eq:rankInequality}, we get:
\begin{equation}\label{eq:mainInequality}
\begin{aligned}
    &2 \big[r_{G'}\big(q(k+n+1)(k+1)\cdot v\big)-\deg(E_2)\big]\geq \\ 
    &\qquad \alpha(\alpha+3)-(\alpha-k)(\alpha-k+1) -l(l+3)+2\gamma+r_0(r_0+1)=:\delta
\end{aligned}
\end{equation}
Set $\beta:=q(k+1)\geq(k+1)$. Now,
\[
l=r_0+q(k+1)+s=r_0+\beta+s \quad\text{and}\quad \alpha(k+1)=\beta(k+n+1).
\]

A direct computation gives
\begin{equation*}
\begin{aligned}
    \alpha(\alpha+3)-(\alpha-k)(\alpha-k+1)%&=(A^2+3A)-\big(A^2+(2k-1)A+(k^2-k)\big) \\
    &= 2(k+1)\alpha-k^2+k=2\beta(k+n+1)-k^2+k,
\end{aligned}
\end{equation*}
and
\begin{equation*}
    2\gamma=(r_0+s)^2+3(r_0+s)-r_0^2-3r_0=2r_0s+s^2+3s.
\end{equation*}
Next,
\begin{equation*}
\begin{aligned}
    l(l+3)&=l^2+3l=(r_0+\beta+s)^2+3(r_0+\beta+s) \\
    &=r_0^2+\beta^2+s^2+2r_0s+2r_0\beta+2s\beta+3r_0+3\beta+3s
\end{aligned}
\end{equation*}
Substituting the results above gives
\begin{equation*}
\begin{aligned}
    \delta&=\big(2\beta(k+n+1)-k^2+k\big)-\big(r_0^2+\beta^2+s^2+2r_0s+2r_0\beta+2s\beta+3r_0+3\beta+3s\big) \\
    &\qquad +(2r_0s+s^2+3s)+(r_0^2+r_0)  \\
    &=-\beta^2+(2k+2n-2r_0-2s-1)\beta-k^2+k-2r_0
\end{aligned}
\end{equation*}
Since, $2r_0=2l-2\beta-2s\leq k+n-2-2\beta-2s$, we get
\begin{equation*}
\begin{aligned}
    \delta&=-\beta^2+(2k+2n-2r_0-2s-1)\beta-k^2+k-2r_0 \\
    &\geq -\beta^2+\big(2k+2n-(k+n-2-2\beta-2s)-2s-1\big)\beta-k^2+k-(k+n-2-2\beta-2s) \\
    %&=-B^2+(k+n+2B+3)B-k^2-n+2s+2 \\
    &=\beta^2+(k+n+3)\beta-k^2-n+2s+2 \\
    &\geq (k+1)^2+(k+n+3)(k+1)-k^2-n+2s+2 \\
    %&=(k^2+2k+1)+(k^2+k+kn+n+3k+3)-k^2-n+2s+2 \\
    &=k^2+6k+kn+2s+6>0.
\end{aligned}
\end{equation*}
Finally, from \cref{eq:mainInequality} and the computation above, we obtain
\begin{equation*}
    2 \big[r_{G'}\big(q(k+n+1)(k+1)\cdot v\big)-\deg(E_2)\big]\geq \delta>0,
\end{equation*}
which proves the claim.
\end{proof}

Now, by \cref{eq:rD_2Equivalence}, we have:
\begin{equation*}
\begin{aligned}
    l\cdot D_2-E&\sim \big((r_0+s)\cdot D_2+q(k+n+1)(k+1)\cdot v\big)-E|_{K_n}-E_1-E_2 \\
    &= \big((r_0+s)\cdot D_2-E_{K_n}-E_1\big)+\big(q(k+n+1)(k+1)\cdot v-E_2\big).
\end{aligned}
\end{equation*}
By \cref{claim:1stEffectiveDivisor,claim:2ndEffectiveDivisor}, both terms on the right-hand side are equivalent to effective divisors, concluding Case (E). 

\medskip

Therefore, for every integer $-2\leq l\leq \frac{k+n-2}{2}$, we have
\begin{equation}\label{eq:lD_2Inequality}
    r_{G}\big(l(k+n+1)\cdot u\big)=r_{G}(l\cdot D_2) \geq  \frac{l(l+3)}{2}.
\end{equation}

To extend this inequality to all $\frac{k+n-2}{2}<l\leq k+n$, we take a different approach.
First, note that by Part (2) of the lemma, 
\[
(k+n-2)\cdot D_2\sim K_G \qquad\text{and}\qquad K_G-l\cdot D_2\sim (k+n-l-2)\cdot D_2.
\]
Now fix some $u\in V(G')$ and integer $\frac{k+n-2}{2}<l\leq k+n$. We have that
\[
-2\leq k+n-l-2<\frac{k+n-2}{2}.
\]
From \cref{eq:lD_2Inequality} and the Riemann--Roch theorem, after some computation, it now follows that
\begin{equation*}
\begin{aligned}
    r_{G}\big(l(k+&n+1)\cdot u\big)=r_{G}(l\cdot D_2)=\deg(l\cdot D_2)-g_G+1+r_{G}(K_G-l\cdot D_2) \\
    &= l(k+n+1)-\binom{k+n}{2}+1+r_{G}\big((k+n-l-2)\cdot D_2\big) \\
    &\geq \frac{2l(k+n)+2l}{2}-\frac{(k+n)^2-(k+n)}{2}+1+\frac{(k+n-l-2)(k+n-l+1)}{2} \\
    %&=\frac{2l(k+n)+2l-(k+n)^2+(k+n)+2}{2} \\
    %&\qquad+\frac{(k+n)^2-(k+n)l+(k+n)-l(k+n)+l^2-l-2(k+n)+2l-2}{2} \\
    &%=\frac{l^2+3l}{2}
    =\frac{l(l+3)}{2}.
\end{aligned}
\end{equation*}
This establishes the inequality for all $-2\leq l\leq k+n$, and proves the lemma.

\subsection{Bounding gonalities from below}
This subsection is dedicated to proving the following lemma. 

\begin{lemma}\label{lemma:lowerBoundMainLemma}
    Let $G'$ be a graph with genus $\binom{k}{2}$ for some $k\geq1$. Assume that for every $1\leq l\leq k-2$,
    \[
    \gon_{\frac{l(l+1)}{2}}(G')\geq lk,
    \]
    and for every integer $r\geq1$, the divisor
    $\gon_r(G')\cdot v$
    realises the $r$-th gonality of $G'$. Let $G$ be the graph obtained by attaching a complete graph $K_n$ to $G'$ along $k+1$ edges between $v\in V(G')$ and each vertex in $K_n$ (see \cref{fig:G} for an illustration). Then, for every $1\leq l\leq k+n-2$,
    \[
    \gon_{\frac{l(l+1)}{2}}(G)\geq l(k+n).
    \]
\end{lemma}

The lemma  builds on the results of Cools and Panizzut \cite{CoolsPanizzut_gonalityComplete}, who compute lower bounds for the gonalities of complete graphs. We summarise their results in the following proposition.
\begin{proposition}[Lower bound for the gonality sequence of complete graphs, \cite{CoolsPanizzut_gonalityComplete}]\label{proposition:lowerBoundForGonalitySequenceOfCompleteGraphs}
    Let $K_n$ be a complete graph with vertices $v_1,v_2,\ldots, v_n$. For any $1\leq l\leq n-3$, let $D$ be a $v_1$-reduced divisor with
    \[
    \deg(D)= l(n-1)-1
    \]
    satisfying
    \[
    D(v_2)\leq D(v_3)\leq\cdots\leq D(v_n).
    \]
    Then, there exists an integer $d\geq 0$ such that firing  $d$ times from $v_1$ produces a divisor $D'$ given by
    \[
    D'(v_i) =
    \begin{cases}
        D(v_1) -d(n-1), & \text{if } i = 1, \\
        D(v_i)+d, & \text{if }2 \le i \le n,
    \end{cases}
    \]
    for which the divisor $E$ defined by
    \[
    E(v_i) =
    \begin{cases}
        \max\big(0,D'(v_1) + 1\big) & \text{if } i = 1, \\
        \max\bigl(0,\, D'(v_i)-i+2\bigr), & \text{if }2 \leq i \leq n,
    \end{cases}
    \]
    has  $\deg(E)\leq\frac{l(l+1)}{2}$. Moreover, $D'-E$ is $v_1$-reduced, has rank $-1$ and satisfies
    \[
    (D'-E)(v_2)\leq(D'-E)(v_3)\leq \cdots\leq(D'-E)(v_n).
    \]
\end{proposition}

For more details, see \cite[Section 4, Remark 12]{CoolsPanizzut_gonalityComplete}.
\begin{lemma}[Lemma 5, \cite{CoolsPanizzut_gonalityComplete}]\label{lemma:reducedDivisorOnCompleteGraph}
    Let $D$ be a divisor on the complete graph $K_n$. Fix a vertex $v_1$, and order the remaining vertices $v_2,v_3,\ldots, v_n$ so that,
    \[
    D(v_2)\leq D(v_3)\leq\ldots\leq D(v_n).
    \]
    Then $D$ is $v_1$-reduced if and only if 
    \[
    0\leq D(v_i)\leq i-2, \qquad\text{for every $2\leq i\leq n$.}
    \]
    
\end{lemma}

We now prove \cref{lemma:lowerBoundMainLemma}. Fix $1\leq l\leq k+n-2$ and let $D$ be an effective divisor on $G$ with $\deg(D)\leq l(k+n)-1$. We aim to show that the rank of $D$ is smaller than $\frac{l(l+1)}{2}$. Since rank is invariant under linear equivalence, we may assume that $D$ is $v$-reduced.
\begin{claim}\label{claim:supportedOnK_n}
    Let $D$ be an effective $v$-reduced divisor on $G$. Define
    \[
    D'=\deg(D|_{G'})\cdot v+D|_{K_{n}}.
    \]
    Then,
    \[
    r_G(D)\leq r_G(D').
    \]
    Note that the divisors $D$ and $D'$ both have the same degree.
\end{claim}
\begin{proof}
    Let $E$ be an effective divisor on $G$ with $\deg(E)=r_G(D)$. We will show that $D'-E$ is effective. Write
    \[
    E=E|_{G'}+E|_{K_n}.
    \]
    As $D$ is $v$-reduced,  \cref{lemma:canRestrictToComponent} implies that the divisor $D|_{K_n\cup\{v\}}$ has rank at least $\deg(E)$ as a divisor on $K_n\cup\{v\}$. Therefore $D|_{K_n\cup \{v\}}-E|_{K_n}$ is equivalent to an effective $v$-reduced divisor $F$ on $K_n\cup \{v\}$.

    We claim that the divisors $D|_{K_n\cup \{v\}}-E|_{K_n}$ and $F$ are equivalent as divisors on $G$ as well. Take any sequence of chip-firing moves on $K_n\cup\{v\}$ that transforms the divisor $D|_{K_n\cup \{v\}}-E|_{K_n}$ to $F$. We can extend it to a sequence of moves on $G$ by replacing every firing of $v$ with  set-firing the whole of $V(G')$. As $v$ is a cut vertex, applying the chip-firing moves in this sequence to the divisor $D|_{K_n\cup \{v\}}-E|_{K_n}$ on $G$ has the exact same effect on $K_n\cup \{v\}$, while keeping the configuration of chips at $G'\setminus \{v\}$ unchanged throughout.
    
    Therefore, we have
    \[
    D-E|_{K_n}=D|_{G'\setminus\{v\}}+D|_{K_n\cup\{v\}}-E|_{K_n}\sim D|_{G'\setminus\{v\}}+F=D|_{G'\setminus\{v\}}+ F(v)\cdot v+F|_{K_n}\geq 0,
    \]
    and
    \[
    \begin{aligned}
    D'-E|_{K_n}&=\deg(D|_{G'\setminus\{v\}})\cdot v + D|_{K_n\cup\{v\}}-E|_{K_n} \\
    &\sim\deg(D|_{G'\setminus\{v\}})\cdot v +F
    =\big(\deg(D|_{G'\setminus\{v\}})+F(v) \big)\cdot v+ F|_{K_n}\geq 0.
     \end{aligned}
    \]
    As $E$ has degree equal to $r_G(D)$,
    \[
    r_G\big(D|_{G'\setminus\{v\}}+ F(v)\cdot v+F|_{K_n} \big)=r_G\big(D-E|_{K_n}\big)\geq r_G(D)-\deg(E|_{K_n})=\deg(E|_{G'}).
    \]
    Notice that since both $F$ and $D$ are $v$-reduced, $D|_{G'\setminus\{v\}}+ F(v)\cdot v+F|_{K_n}$ is $v$-reduced as well. So by \cref{lemma:canRestrictToComponent}, 
    \[
    r_{G'}\big(D|_{G'\setminus\{v\}}+ F(v)\cdot v\big)\geq \deg(E|_{G'}).
    \]
    Therefore, 
    \[
    \gon_{\deg(E|_{G'})}(G')\leq\deg\big(D|_{G'\setminus\{v\}}+ F(v)\cdot v\big)=\deg(D|_{G'\setminus\{v\}})+F(v).
    \]
    By assumption, 
    \[
    \deg(E|_{G'})=r_{G'}\big(\gon_{\deg(E|_{G'})}(G')\cdot v\big)\leq r_{G'}\big[\big(\deg\left(D|_{G'\setminus\{v\}}\right)+F(v) \big)\cdot v\big].
    \]
    Now we can apply \cref{lemma:restrictedRank} %\yoav{reference?} 
    to $\big(\deg(D|_{G'\setminus\{v\}})+F(v) \big)\cdot v+ F|_K$, concluding that 
    \[
    \big[\deg\big(D|_{G'\setminus\{v\}}\big)+F(v) \big]\cdot v+ F|_K-E|_{G'}
    \]
    is equivalent to an effective divisor on $G$. Now,
    \[
    D'-E=D'-E|_K-E|_{G'}\sim \big[\deg\big(D|_{G'\setminus\{v\}}\big)+F(v) \big]\cdot v+ F|_K-E|_{G'}
    \]
    is equivalent to an effective divisor as well. As $E$ was an arbitrary effective divisor of degree $r_G(D)$, the claim is follows.
\end{proof}

By \cref{claim:supportedOnK_n}, it is enough to show that the rank of $D'=\deg(D|_{G'})\cdot v+D|_{K_n}$ is smaller than $\frac{l(l+1)}{2}$. 
Let $u_1$ be any vertex in $K_n$. Perform Dhar's burning algorithm to obtain a $u_1$-reduced divisor $D_G$ on $G$ equivalent to $D'$.
As there are no chips on $G'\setminus\{v\}$, during each iteration of the algorithm, either all of $G'$ burns down or none of it. 
So the number of chips on vertices in $G'$ other than $v$ does not change. In particular, $D_G$ is supported on $K_n\cup \{v\}$.
Label the other vertices of $K_n$ by $u_2,u_3,\ldots, u_n$ so that
\begin{equation}\label{eq:orderingChipsOnK_n}
D_G(u_2)\leq D_G(u_3)\leq \cdots\leq D_G(u_n).
\end{equation}
Define 
\[
\alpha=\left\lfloor \frac{D_G(v)}{k+1}\right\rfloor
\]
and let $0\leq \beta\leq k$ be the remainder of the division of $D_G(v)$ by $k+1$. So
\[
D_G(v)=\alpha(k+1)+\beta.
\]
Set
\[
N=\max\big(\{1\}\cup\{2\leq i\leq n \mid D_G(u_i)\leq\alpha\}\big).
\]
Construct a divisor $D_K$ on the complete graph $K_{n+k+1}$ with vertices $v_1,v_2,\ldots, v_{n+k+1}$ by
\[
D_K(v_i)=
\begin{cases}
    D_G(u_i) & \text{if } 1\leq i\leq N, \\
    \alpha  & \text{if } N+1\leq i\leq N+k+1-\beta, \\
    \alpha+1 & \text{if } N+k+2-\beta\leq i\leq N+k+1, \\
    D_G\left(u_{i-(k+1)}\right) & \text{if } N+k+2\leq i\leq n+k+1.
\end{cases}
\]
See \cref{fig:D_G and D_K} for an example.

\begin{figure}
\centering
\begin{tikzpicture}[scale=1.75, node/.style={circle, draw, fill=black, inner sep=1.2pt}]
    % Vertex v
    \node[node, label=left:{$12$}] (v) at (0.3,0) {};

    % Graph G'
    \draw (-0.6,0) circle [radius=1];
    %\node at (-0.6,0) {$G'$};

    % Complete graph K_n
    %\node at (2.3,1.1) {$K_n$};
    \node[node, label=above:{$u_1$ $19$}] (u1) at (1.5,1) {};
    \node[node, label=above right:{$0$}] (u2) at (2.3,0.5) {};
    \node[node, label=below right:{$1$}] (u3) at (2.3,-0.5) {};
    \node[node, label=below:{$7$}] (u4) at (1.5,-1) {};

    % Edges of K_n
    \draw (u1) -- (u2);
    \draw (u1) -- (u3);
    \draw (u1) -- (u4);
    \draw (u2) -- (u3);
    \draw (u2) -- (u4);
    \draw (u3) -- (u4);

    % Edges between v and K_n
    \draw[bend right=0](v) to (u1);
    \draw[bend left=5]  (v) to (u1);
    \draw[bend right=5] (v) to (u1);
    \draw[bend left=10] (v) to (u1);
    \draw[bend right=10](v) to (u1);

    \draw[bend right=0](v) to (u2);
    \draw[bend left=5]  (v) to (u2);
    \draw[bend right=5] (v) to (u2);
    \draw[bend left=10] (v) to (u2);
    \draw[bend right=10](v) to (u2);
    
    \draw[bend right=0](v) to (u3);
    \draw[bend left=5]  (v) to (u3);
    \draw[bend right=5] (v) to (u3);
    \draw[bend left=10] (v) to (u3);
    \draw[bend right=10](v) to (u3);

    \draw[bend right=0](v) to (u4);
    \draw[bend left=5]  (v) to (u4);
    \draw[bend right=5] (v) to (u4);
    \draw[bend left=10] (v) to (u4);
    \draw[bend right=10](v) to (u4);

    % Complete graph K_9
    \begin{scope}[shift={(5,0)}]
        % Radius
        \def\r{1.2}
        % Vertex labels
        \def\labone{$v_1$ $19$}
        \def\labtwo{$0$}
        \def\labthree{$1$}
        \def\labfour{$2$}
        \def\labfive{$2$}
        \def\labsix{$2$}
        \def\labseven{$3$}
        \def\labeight{$3$}
        \def\labnine{$7$}

        % Vertices
        \node[node, label=above:{\labone}]
            (w1) at ({-\r*cos(50 + 40*1)}, {\r*sin(50 + 40*1)}) {};
        \node[node, label=above right:{\labtwo}]
            (w2) at ({-\r*cos(50 + 40*2)}, {\r*sin(50 + 40*2)}) {};
        \node[node, label=right:{\labthree}]
            (w3) at ({-\r*cos(50 + 40*3)}, {\r*sin(50 + 40*3)}) {};
        \node[node, label=below right:{\labfour}]
            (w4) at ({-\r*cos(50 + 40*4)}, {\r*sin(50 + 40*4)}) {};
        \node[node, label=below:{\labfive}]
            (w5) at ({-\r*cos(50 + 40*5)}, {\r*sin(50 + 40*5)}) {};
        \node[node, label=below:{\labsix}]
            (w6) at ({-\r*cos(50 + 40*6)}, {\r*sin(50 + 40*6)}) {};
        \node[node, label=below left:{\labseven}]
            (w7) at ({-\r*cos(50 + 40*7)}, {\r*sin(50 + 40*7)}) {};
        \node[node, label=left:{\labeight}]
            (w8) at ({-\r*cos(50 + 40*8)}, {\r*sin(50 + 40*8)}) {};
        \node[node, label=above left:{\labnine}]
            (w9) at ({-\r*cos(50 + 40*9)}, {\r*sin(50 + 40*9)}) {};

        % Draw edges
        \foreach \i in {1,...,9} {
            \foreach \j in {\i,...,9} {
                \ifnum\i<\j
                    \draw (w\i) -- (w\j);
                \fi
            }
        }

        %\node at (0,1.8) {$K_9$};
    \end{scope}
\end{tikzpicture}

\caption{An example of a $u_1$-reduced divisor $D_G$ on the left and the corresponding divisor $D_K$ on the right. Here, $N=3$ and $k=n=4$. }
\label{fig:D_G and D_K}
\end{figure}

\begin{claim}
    The divisor $D_K$ is $v_1$-reduced on $K_{n+k+1}$, with $\deg(D_K)=\deg(D)=l(k+n)-1$ and
    \[
    D_K(v_2)\leq D_K(v_3)\leq\cdots\leq D_K(v_{n+k+1}).
    \]
\end{claim}
\begin{proof}
    We first verify the degree of $D_K$.
    \begin{equation*}
    \begin{aligned}
        \deg(D_K)&=\sum_{i=1}^ND_K(v_i)+\sum_{i=N+1}^{N+k+1-\beta}D_K(v_i)+\sum_{i=N+k+2-\beta}^{N+k+1}D_K(v_i)+\sum_{i=N+k+2}^{n+k+1}D_K(v_i) \\
        &=\sum_{i=1}^ND_G(u_i)+\sum_{i=N+1}^{N+k+1-\beta}\alpha+\sum_{i=N+k+2-\beta}^{N+k+1}(\alpha+1)+\sum_{i=N+1}^{n}D_G(u_i) \\
        &=\alpha(k+1)+\beta+\sum_{w \in V(K_n)}D_G(w) \\
        &=\deg(D_G|_{G'})+\deg(D_G|_{K_n})=\deg(D_G)=l(k+n)-1.
    \end{aligned}  
    \end{equation*}
    Since $D_K(v_i)=D_G(u_i)$ for every $2\leq i\leq N$, and $D_K(v_i)=D_G(u_{i-(k+1)})$ for every $N+k+2\leq i\leq n+k+1$, from \cref{eq:orderingChipsOnK_n} we have
    \[
    D_K(v_2)\leq D_K(v_3)\leq\cdots\leq D_K(v_N)\qquad \text{and}\qquad D_K(v_{N+k+2})\leq D_K(v_{N+k+3})\leq\cdots\leq D_K(v_{n+k+1}).
    \]
    By construction,
    \[
    D_K(v_{N+1})\leq D_K(v_{N+2})\leq\ldots\leq D_K(v_{N+k+1}).
    \]
    and by definition of $N$,
    \[
    D_K(v_{N+1})=\alpha\geq D_G(u_{N})=D_K(v_{N}), \qquad \text{if $N\geq 2$},
    \]
    and
    \[
    D_K(v_{N+k+1})\leq\alpha+1\leq D_G(u_{N+1})=D_K(v_{N+k+2}), \qquad \text{if $N\leq n-1$.}
    \]
    These inequalities give
    \[
    D_K(v_2)\leq D_K(v_3)\leq \ldots\leq D(v_{n+k+1}).
    \]
    By \cref{lemma:reducedDivisorOnCompleteGraph}, it remains to show that $0\leq D(v_i)\leq i-2$ for all $2\leq i\leq n+k+1$.
    Assume for a contradiction that there exists $2\leq j\leq n+k+1$ such that $D_K(v_j)\geq j-1$. We consider 3 cases.

    \casetitle[$2\leq j\leq N$]{A}

    \medskip
    \noindent
    Let $A=V\left(G'\right)\cup\{u_j,u_{j+1},\ldots, u_{n}\}$. Then,
        \[
        \operatorname{outdeg}_A(v)=(j-1)(k+1)\leq D_K(v_{k+1})(k+1)\leq \alpha(k+1)\leq D_G(v),
        \]
        and
        \[
        \operatorname{outdeg}_{A}(u_i)=j-1\leq D_K(v_j)=D_G(u_j)\leq D_G(u_i),\qquad \text{for all $j\leq i\leq n$.}
        \]
        Therefore, set-firing $A$ with respect to $D_G$ doesn't create any debt on $G\setminus\{u_1\}$, contradicting that $D_G$ is $u_1$-reduced.

    \casetitle[$N+1\leq j\leq N+k+1$]{B}

    \medskip
    \noindent
    Let $A=V(G')\cup\left\{u_{N+1},u_{N+2}\ldots u_{n}\right\}$.
        Notice that if $j=N+1$, then $\alpha =D_K(v_{N+1})\geq N$, and if $N+2\leq j\leq N+k+1$,
        \[
        \alpha=D_K(v_{N+1})\geq D_K(v_{j})-1\geq j-2\geq N.
        \]
        Now, 
        \[
        \operatorname{outdeg}_A(v)=N(k+1) \leq \alpha(k+1)\leq D_G(v),
        \]
        and
        \[
        \begin{aligned}
        \operatorname{outdeg}_{A}(u_i)&=N\leq j-1\leq D_K(v_j) \\
        &\leq D_K(v_{N+k+2})= D_G(u_{N+1})
        \leq D(u_{i}),\qquad \text{for all $N+1\leq i\leq n$.}
        \end{aligned}
        \]
        Again, set-firing $A$ with respect to $D_G$ creates no debt on $G\setminus\{u_1\}$, contradicting that $D_G$ is $u_1$-reduced.

    \casetitle[$N+k+2\leq j\leq n+k+1$]{C}

    \medskip
    \noindent
    Let $A=\left\{u_{j-(k+1)},u_{j-(k+1)+1},\ldots, u_n\right\}$.
        \[
        \operatorname{outdeg}_{A}(u_i)=j-1\leq D_K(v_j)=D_G(u_j)\leq D_G(u_i),\qquad \text{for all $j-(k+1)\leq i\leq n$.}
        \]
        Set-firing $A$ with respect to $D_G$ doesn't create any debt on $G\setminus\{u_1\}$, again contradicting that $D_G$ is $u_1$-reduced.

    \medskip
    We obtain a contradiction in all cases, showing no such $j$ exists. By \cref{lemma:reducedDivisorOnCompleteGraph}, $D_K$ is $v_1$-reduced.
\end{proof}

We now apply \cref{proposition:lowerBoundForGonalitySequenceOfCompleteGraphs} to the divisor $D_K$. This yields a non-negative integer $d$ and divisors $D_K'$ and $E_K$  as described in the proposition. Explicitly,
\[
D'_K(v_i) =
\begin{cases}
    D_K(v_i)-d(k+n) & \text{if } i=1, \\
    D_K(v_i)+d & \text{if } 2\leq i\leq N, \\
    \alpha+d  & \text{if } N+1\leq i\leq N+k-\beta+1, \\
    \alpha+d+1 & \text{if } N+k-\beta+2\leq i\leq N+k+1, \\
    D_K(v_i)+d & \text{if } N+k+2\leq i\leq n+k+1,
\end{cases}
\]
and
\[
E_K(v_i) =
\begin{cases}
    \max\big(0,D_K'(v_1) + 1\big) & \text{if } i = 1, \\
    \max\bigl(0,\, D_K'(v_i)-i+2\bigr) & \text{if }2 \leq i \leq n+k+1,
\end{cases}
\]
with $\deg(E_K)\leq \frac{l(l+1)}{2}$. Let 
\[
r:=E_K(v_{N+1})=\max\bigl(0,D'(v_{N+1})-N+1\bigr)=\max(0,\, \alpha+d-N+1)\geq0.
\]
Note that if $r>0$, then $N=\alpha+d-r+1$.
\begin{claim}
Let $H$ be a subgraph of $K_{n+k+1}$ consisting of vertices $v_{N+1}, v_{N+2}\ldots v_{N+k+1}$ and all edges between them.
    Then
    \[
    \deg(E_K|_{H})=
    \begin{cases}
        \frac{r(r+1)}{2} & \text{if } 0\leq r\leq k-\beta,\\
        \frac{r(r+1)}{2}+\beta-r+k  & \text{if } k-\beta+1\leq r\leq k-1, \\
        \frac{k(k+1)}{2}+(k+1)(r-k)+\beta  & \text{if } k\leq r.
    \end{cases}
    \]
\end{claim}
\begin{proof}
First notice that if $E_K(v_i)=0$ for some $N+1\leq i \leq N+k+1$, then for any $i+1\leq j \leq N+k+1$ we have
\[
\begin{aligned}
    0\geq D'_K(v_{i})-i+2\geq D'_K(v_j)-j+2\geq E_K(v_j)\geq 0,
\end{aligned}
\]
hence $E_K(v_j)=0$. We consider the following four cases.

\casetitle[$r=0$]{A}

    \medskip
    \noindent
    Then $E_K(v_{N+1})=0$, and so $E_K(v_j)=0$ for all $N+2\leq j \leq N+k+1$. Therefore, $\deg(E_K|_H)=0$.

\casetitle[$1\leq r\leq k-\beta$]{B}

    \medskip
    \noindent
    For each $N+1\leq i\leq N+r+1\leq N+k-\beta+1$,
    \[
    \begin{aligned}
        E_K(v_i)&=\max\big(0,D'_K(v_i)-i+2\big)=\max(0,\alpha+d-i+2) \\
        &=\max(0,N+r+1-i)=N+r+1-i.
    \end{aligned}
    \]
    In particular, $E(v_{N+r+1})=0$. Hence,
    \[
    \deg(E_K|_{H})=\sum_{i=N+1}^{N+r+1}E_K(v_i)=\sum_{i=N+1}^{N+r+1} (N+r+1-i)=\sum_{j=0}^{r}j=\frac{r(r+1)}{2}.
    \]

    \casetitle[$k-\beta+1\leq r\leq k-1$]{C}

    \medskip
    \noindent
    For every $N+1\leq i\leq N+k-\beta+1< N+r+1$,
    \[
    \begin{aligned}
        E_K(v_i)&=\max\big(0,D'_K(v_i)-i+2\big)=\max(0,\alpha+d-i+2) \\
        &=\max(0,N+r+1-i)=N+r+1-i,
    \end{aligned}
    \]
    and for $N+k-\beta+2\leq i\leq N+r+2\leq N+k+1$,
    \[
    \begin{aligned}
        E_K(v_i)&=\max\big(0,D'_K(v_i)-i+2\big)=\max(0,\alpha+d+1-i+2) \\
        &=\max(0,N+r+2-i)=N+r+2-i.
    \end{aligned}
    \]
    In particular, $E_K(v_{N+r+2})=0$. Then
    \[
    \begin{aligned}
    \deg(E_K|_{H})&=\sum_{i=N+1}^{N+r+2}E_K(v_i)= \sum_{i=N+1}^{N+k-\beta+1} E_K(v_i)+ \sum_{i=N+k-\beta+2}^{N+r+2} E_K(v_i)\\ &=\sum_{i=N+1}^{N+k-\beta+1} (N+r+1-i)+ \sum_{i=N+k-\beta+2}^{N+r+2} (N+r+2-i)\\
    &=\sum_{j=r-k+\beta}^{r}j+\sum_{j=0}^{r-k+\beta-1}(j+1)=\sum_{j=0}^{r}j+(r-k+\beta)=\frac{r(r+1)}{2}+r-k+\beta,
    \end{aligned}
    \]
    as claimed.

    \casetitle[$k\leq r$]{D}

    \medskip
    \noindent
    For all $N+1\leq i\leq N+k-\beta+1\leq N+r+1$,
    \[
    \begin{aligned}
        E_K(v_i)&=\max\big(0,D'_K(v_i)-i+2\big)=\max(0,\alpha+d-i+2) \\
        &=\max(0,N+r+1-i)=N+r+1-i,
    \end{aligned}
    \]
    and for every $N+k-\beta+2\leq i\leq N+k+1< N+r+2$,
    \[
    \begin{aligned}
        E_K(v_i)&=\max\big(0,D'_K(v_i)-i+2\big)=\max(0,\alpha+d+1-i+2) \\
        &=\max(0,N+r+2-i)=N+r+2-i.
    \end{aligned}
    \]
    Therefore,
    \[
    \begin{aligned}
    \deg(E_K|_{H})&=\sum_{i=N+1}^{N+r+2}E_K(v_i)= \sum_{i=N+1}^{N+k-\beta+1} E_K(v_i)+ \sum_{i=N+k-\beta+2}^{N+k+1} E_K(v_i)\\ 
    &=\sum_{i=N+1}^{N+k-\beta+1} (N+r+1-i)+ \sum_{i=N+k-\beta+2}^{N+k+1} (N+r+2-i)\\
    &=\sum_{j=r-k+\beta}^{r}j+\sum_{j=r-k}^{r-k+\beta-1}(j+1)=\sum_{j=\beta}^{k}(j+r-k)+\sum_{j=0}^{\beta-1}(j+r-k+1) \\
    &=\sum_{j=0}^{k}j+(k-\beta+1)(r-k)+\beta(r-k+1)=\frac{k(k+1)}{2}+(k+1)(r-k)+\beta,
    \end{aligned}
    \]
    which completes the proof.
\end{proof}

We now return to the graph $G$ and divisor $D_G$. Let $D_G^\prime$ be the divisor obtained by firing $d$ times at the vertex $u_1$ with respect to $D_G$. Then
\[
D'_G(w) =
\begin{cases}
    D_G(u_1)-d(k+n) & \text{if } w=u_1, \\
    D_G(w)+d & \text{if } w\in\{u_2,\ldots, u_N\}, \\
    (\alpha+d)(k+1)+\beta  & \text{if } w=v, \\
    0 & \text{if } w\in V(G')\setminus\{v\}, \\
    D_G(w)+d & \text{if } w\in\{u_{N+1},\ldots, u_n\}.
\end{cases}
\]
Notice that for $1\leq i\leq N$, 
\[
D'_G(u_i)=D'_K(v_i)
\]
and for $N+1\leq i\leq n$, 
\[
D'_G(u_i)=D'_K(v_{i+k+1}).
\]

\begin{claim}\label{claim:divisorsE_rAndD_1}
    There is an effective divisor $E_r$ on $G'$ such that  $\deg(E_r)=\deg(E_K|_H)$ and $D'_G|_{G'}-E_r$ is equivalent to a $v$-reduced divisor $D_1$ on $G'$ satisfying
    \[
    D_1(v)<N(k+1).
    \]
\end{claim}
\begin{proof}
    First, note that
    \[
    D'_G|_{G'}=\big((\alpha+d)(k+1)+\beta\big)\cdot v.
    \]
    We regard this as a divisor on $G'$ throughout the proof.
    If $l=k-1$ or $k$, then 
    \[
    \frac{l(l+1)}{2}\geq \frac{(k-1)k}{2}=g_{G'}.
    \]
    So by \cref{corollary:corollariesOfRiemann-Roch} we have,
    \[
    \begin{aligned}
        \gon_{\frac{(k-1)k}{2}}(G')=\gon_{g_{G'}}(G')=2g_{G'}=(k-1)k, \\
        \gon_{\frac{k(k+1)}{2}}
        (G')=\gon_{g_{G'}+k}(G')=2g_{G'}+k=k^2.
    \end{aligned}
    \]
    Similarly as in the proof of the previous claim, we consider four cases.

    \casetitle[$r=0$]{A}

    \medskip
    \noindent
    Since $r=\max(0,\, \alpha+d-N+1)$, we have $\alpha+d-N+1\leq 0$.
    Hence,
    \[
    D'_G|_{G'}(v)=(\alpha+d)(k+1)+\beta<(\alpha+d+1)(k+1)=N(k+1).
    \]
    Taking $D_1=D'_G|_{G'}$ and $E_r=0$ concludes this case.

    \casetitle[$1\leq r\leq k-\beta$]{B}

    \medskip
    \noindent
    Since $\gon_{\frac{r(r+1)}{2}}(G')\geq rk$,
    there exists an effective divisor $E_r$ on $G'$ with degree
    \[
    \deg(E_r)=\frac{r(r+1)}{2}=\deg(E_K|_H)
    \]
    such that the divisor $(rk-1)\cdot v-E_r$ is unwinnable. Let $D_0$ be the $v$-reduced divisor equivalent to $(rk-1)\cdot v-E_r$, so $D_0(v)\leq-1$. Define $D_1$ as,
    \begin{equation*}
    \begin{aligned}
        D'_G|_{G'}-E_r&=\big((\alpha+d)(k+1)+\beta\big)\cdot v-E_r \\
        &=\big((\alpha+d)(k+1)+\beta-rk+1\big)\cdot v +(rk-1)\cdot v -E_r \\
        &\sim \big((\alpha+d)(k+1)+\beta-rk+1\big)\cdot v +D_0 =:D_1.
    \end{aligned}
    \end{equation*}
    Then $D_1$ is $v$-reduced and satisfies
    \[
    D_1(v)\leq (\alpha+d)(k+1)+\beta -rk.
    \]
    Now,
    \begin{equation*}
    \begin{aligned}
        N(k+1)-D_1(v)&\geq (\alpha+d-r+1)(k+1)-\big((\alpha+d)(k+1)+\beta -rk\big) \\
        &= k-r-\beta +1\geq 1>0.
    \end{aligned}
    \end{equation*}

    \casetitle[$k-\beta+1\leq r\leq k-1$]{C}

    \medskip
    \noindent
    We have
    \[
    \gon_{\frac{r(r+1)}{2}+\beta-r+k}(G')\geq \gon_{\frac{r(r+1)}{2}}(G') +\beta-r+k\geq rk+\beta -r +k.
    \]
    Therefore, there exists an effective divisor $E_r$ with degree 
    \[
    \deg(E_r)=\frac{r(r+1)}{2}+\beta-r+k=\deg(E_K|_H)
    \]
    such that the divisor
    \[
    (rk+\beta-r+k-1)\cdot v-E_r
    \]
    is unwinnable. That is, it is equivalent to a $v$-reduced divisor $D_0$ with $D_0(v)\leq -1$. Then,
    \[
    \begin{aligned}
        D'_G|_{G'}-E_r&=\big((\alpha+d)(k+1)+\beta\big)\cdot v-E_r \\
        &=\big((\alpha+d)(k+1)-rk+r-k+1\big)\cdot v +(rk+\beta-r+k-1)\cdot v-E_r \\
        &\sim \big((\alpha+d)(k+1)-rk+r-k+1\big)\cdot v +D_0 =:D_1.
    \end{aligned}
    \]
    $D_1$ is also $v$-reduced with $D_1(v)\leq (\alpha+d)(k+1)-rk+r-k$. So
    \[
    \begin{aligned}
        N(k+1)-D_1(v)&\geq (\alpha +d-r+1)(k+1)-\big((\alpha+d)(k+1)-rk+r-k\big) \\
        &=2(k-r)+1>0.
    \end{aligned}
    \]

    \casetitle[$k\leq r$]{D}

    \medskip
    \noindent
    Since
    \[
    \gon_{\frac{k(k+1)}{2}+(k+1)(r-k)+\beta}(G')\geq \gon_{\frac{k(k+1)}{2}}(G') +(k+1)(r-k)+\beta\geq k^2+(k+1)(r-k)+\beta,
    \]
    there exists an effective divisor $E_r$ with 
    \[
    \deg(E_r)=\frac{k(k+1)}{2}+(k+1)(r-k)+\beta=\deg(E_K|_H),
    \]
    such that the divisor
    \[
    \big(k^2+(k+1)(r-k)+\beta-1\big)\cdot v-E_r
    \]
    is unwinnable. Let $D_0$ be the corresponding $v$-reduced divisor with $D_0(v)\leq -1$. Then
    \[
    \begin{aligned}
        D'_G|_{G'}-E_r&=\big((\alpha+d)(k+1)+\beta\big)\cdot v-E_r \\
        &=\big((\alpha+d)(k+1)-k^2-(k+1)(r-k)+1\big)\cdot v \\
        &\quad +\big(k^2+(k+1)(r-k)+\beta-1\big)\cdot v-E_r \\
        &\sim \big((\alpha+d)(k+1)-k^2-(k+1)(r-k)+1\big)\cdot v  +D_0 =:D_1.
    \end{aligned}
    \]
    $D_1$ is again $v$-reduced with $D_1(v)\leq (\alpha+d)(k+1)-rk+r-k$. Now,
    \[
    \begin{aligned}
        N(k+1)-D_1(v)&\geq (\alpha +d-r+1)(k+1)-\big((\alpha+d)(k+1)-k^2-(k+1)(r-k)\big) \\
        &=(1-k)(k+1)+k^2=1>0.
    \end{aligned}
    \]

This completes the proof of the claim.
\end{proof}

Define a divisor $E_G$ on $G$ by
\[
E_G(w)=
\begin{cases}
    E_r(w) & \text{if } w\in V(G'), \\
    E_K(v_i) & \text{if $w=u_i$, for $1\leq i\leq N$}, \\
    E_K(v_{i+k+1}) & \text{if $w=u_i$, for $N+1\leq i\leq n$.}
\end{cases}
\]
Then $E_G|_{G'}=E_r$ and $E_G=E_r+E_G|_{K_n}$. A direct computation gives 
\[
\begin{aligned}
\deg(E_G)&=\sum_{i=1}^NE_K(v_i)+\deg(E_r)+\sum_{i=N+k+2}^{n+k+1}E_K(v_i) \\
&=\sum_{i=1}^NE_K(v_i)+\sum_{i=N+1}^{N+k+1}E_K(v_i)+\sum_{i=N+k+2}^{n+k+1}E_K(v_i) \\
&=\deg(E_K)\leq \frac{l(l+1)}{2}.
\end{aligned}
\]
By \cref{claim:divisorsE_rAndD_1}, we have $D'_G|_{K_n}-E_r\sim D_1$ as divisors on $G'$. Take any chip-firing sequence that transforms $D'_G|_{K_n}-E_r$ into $D_1$. We extend this sequence of chip-firings to $G$ by replacing each firing of $v$ with set-firing all of $V(K_n)\cup \{v\}$. Since $v$ is a cut vertex, applying this sequence to $D'_G|_{K_n}-E_r$ as a divisor on $G$, has the same effect on $G'$ and doesn't cause any changes on the configuration of chips on $K_n$. Therefore $D'_G|_{K_n}-E_r$ and $D_1$ are equivalent as divisors on $G$. We now have
\[
D_G-E_G\sim D'_G-E_G=D'_G|_{G'}+D'_G|_{K_n}-E_r-E_G|_{K_n}\sim D_1+ (D'_G-E_G)|_{K_n}=:F'.
\]
Moreover,
\[
F'(u_i)=(D'_G-E_G)|_{K_n}(u_i)=
\begin{cases}
    (D'_K-E_K)(v_i) & \text{if $1\leq i\leq N$}, \\
    (D'_K-E_K)(v_{i+k+1}) & \text{if $N+1\leq i\leq n$.}
\end{cases}
\]
\begin{claim}
    Starting a fire from $u_1$ with respect to $F'$ burns the entire graph $G$. In particular, $D_G-E_G\sim F'$ is unwinnable and 
    \[
    r_G(D_G)<\deg (E_G) \leq\frac{l(l+1)}{2}.
    \]
\end{claim}
\begin{proof}
    Recall that the divisor $D'_K-E_K$ is $v_1$-reduced, with 
    \[
    F'(u_1)=(D'_K-E_K)(v_1)\leq -1,
    \]
    and that
    \[
    (D'_K-E_K)(v_2)\leq (D'_K-E_K)(v_3)\leq \cdots\leq (D'_K-E_K)(v_{n+k+1}).
    \]
    By \cref{lemma:reducedDivisorOnCompleteGraph}, for $2\leq i\leq N$,
    \begin{equation}\label{eq:boundingF'(u_i)1}
        F'(u_i)=(D'_K-E_K)(v_i)\leq i-2,
    \end{equation}
    and for $N+1\leq i\leq n$,
    \begin{equation}\label{eq:boundingF'(u_i)2}
        F'(u_i)=(D'_K-E_K)(v_{i+k+1})\leq i+k-1.
    \end{equation}
    By \cref{claim:divisorsE_rAndD_1},
    \begin{equation}\label{eq:boundingF'(v)}
        F'(v)=D_1(v)<N(k+1),
    \end{equation}
    and as $F'|_{G'}=D_1$, a fire at $v$ with respect to $F'$ burns all of $G'$.

    Now start a fire at $u_1$ with respect to $F'$. From \cref{eq:boundingF'(u_i)1}, fire spread s successively through $u_2,u_3,\ldots, u_N$.
    By \cref{eq:boundingF'(v)}, the vertex $v$ will catches fire, burning down all of $G'$. Finally, from \cref{eq:boundingF'(u_i)2}, the fire continues spreading through $u_{N+1},u_{N+2},\ldots, u_{n}$. As $F'(u_1)\leq-1$, $F'$ is unwinnable, completing the proof.
\end{proof}

Therefore, for any divisor $D$ on $G$ with degree $l(k+n)-1$, we have constructed a divisor $D_G$ such that
\[
r_G(D)\leq r_G(D_G)<\frac{l(l+1)}{2}.
\]
This shows that $\gon_{\frac{l(l+1)}{2}}(G)\geq l(k+n)$, and concludes the proof of the \cref{lemma:lowerBoundMainLemma}.

\section{Proof of the main theorems}\label{sec:proofOfMainTheorem} 

In this section, we use the tools developed throughout the paper to prove the two main results: \cref{mainTheorem:alternativeDefinitonForQBGraphs} and \cref{mainTheorem:propertiesOfQBGraphs}. First, we record the following claim about equivalent divisors on quasi-banana graphs.

\begin{claim}\label{claim:chipFiringOnQuasiBananaGraphs}
    Let $QB_m(N)$ be a quasi-banana graph
    with notation adapted from \cref{def:quasiBananaGraph}. 
    Let $1\leq L\leq K\leq m$. Then, for any $1\leq j\leq q_L$,
    \begin{equation*}
        u_{0,1}+\sum_{w\in C_{0,K}} w \;\sim\; \left(\sum_{i=1}^{L} q_i+2\right)\cdot u_{L,j} + \sum_{w\in C_{L+1,K}} w,
    \end{equation*}
    where
    \begin{equation*}
        C_{A,B}=\bigcup_{i=A}^B C_i.
    \end{equation*}
\end{claim}
\begin{proof}
    Start with the divisor
    \[
    u_{0,1}+\sum_{w\in C_{1,K}} w.
    \]
    Set-fire the sets
    \[
    C_{0,1}\setminus\{u_{1,1}\}, C_{0,2}\setminus\{u_{2,1}\},\ldots, C_{0,L-1}\setminus\{u_{L-1,1}\}
    \] 
    to obtain
    \[
    \left(\sum_{i=1}^{L-1} q_i+2\right)\cdot u_{L-1,1} + \sum_{w\in C_{L,K}} w.
    \]
    Now, if $j=1$, set-firing $C_{0,L}\setminus\{u_{L,1}\}$ produces the required divisor. If $j\geq 2$, set-firing $V(QB_m(N))\setminus\{u_{L,j}\}$ gives the required divisor.
\end{proof}

\subsection{Proof of \cref{mainTheorem:propertiesOfQBGraphs}}\label{sec:proofOfTheoremB} We aim to show that every quasi-banana graph $QB_m(N)$ satisfies the conclusions listed in the theorem. In addition to that, we show that the divisor 
\[
u_{0,1}+\sum_{w\in V\left(QB_m(N)\right)}w
\]
realises the second gonality. 

We proceed by induction on $m$. For the base case $m=0$, there is only one quasi-banana graph $QB_0$ consisting of a single vertex and no edges. The genus of $QB_0$ is $0$,
and by \cref{corollary:corollariesOfRiemann-Roch}, its gonality sequence is $1,2,3,\ldots$. The unique degree $2$ divisor $2u_{0,1}$ realises the second gonality, and
\[
\big(\gon_1(QB_0)-2\big)\cdot2u_{0,1} = -2u_{0,1} = K_{QB_0}
\]
has rank $-1$. Therefore, the theorem holds for $m=0$.

Now assume that $m\geq1$ and that for every quasi-banana graph $QB_{m-1}(N)$, the divisor 
\[
u_{0,1}+\sum_{w\in V\left(QB_{m-1}(N)\right)}w
\]
realises the second gonality, and that all conclusions of the theorem hold. Let 
\[
G=QB_m(q_1,q_w,\ldots, q_m)
\qquad \text{and}\qquad
G'=QB_{m-1}(q_1,q_w,\ldots, q_{m-1}),
\]
and denote $v=u_{m-1,1}$. Define
\[
k_{m-1}=\sum_{i=1}^{m-1} q_i+1,\qquad k_m=\sum_{i=1}^mq_i+1.
\]
The graph $G$ is created from $G'$ by attaching a complete graph $K$ with vertices $u_{m,1},u_{m,2},\ldots, u_{m,q_m}$ along $k_{m-1}+1$ edges between $v$ and each vertex in $K$. We can therefore think of $G'$ as a subgraph of $G$. We define
\[
D_{m-1}=u_{0,1}+\sum_{w\in V(G')} w, \qquad\text{and}\qquad D_m =u_{0,1}+\sum_{w\in V(G)}w.
\]
By the inductive hypothesis, $G'$ is gonality-tight, so
\[
\gon_1(G')=\deg(D_{m-1})-1=|V(G')|=\sum_{i=1}^{m-1} q_i+1=k_{m-1}.
\]
By \cref{claim:chipFiringOnQuasiBananaGraphs},
\[
D_{m-1}\sim \left(\sum_{i=1}^{m-1} q_i+2\right)\cdot u_{m-1,1} = (k_{m-1}+1)\cdot v
\]
and
\[
D_m\sim \left(\sum_{i=1}^{m} q_i+2\right)\cdot u_{m,1} = (k_{m}+1)\cdot u_{m,1},
\]
where $D_{m-1}$ and $D_m$ are divisor on $G'$ and $G$ respectively.
Furthermore, by the inductive hypothesis, 
\[
g_ {G'}=\binom{k_{m-1}}{2},
\]
and for all $1\leq l\leq k_{m-1}$,
\[
r_{G'}\big(l(k_{m-1}+1)\cdot v\big)=r_{G'}(l\cdot D_{m-1})=\frac{l(l+3)}{2}.
\]
Now we can apply \cref{lemma:upperBoundMainLemma} to $G'$ and $G$, with $k=k_{m-1}$ and $n=q_m$. This gives
\begin{equation}\label{eq:genusOfQuasiBananaGraph}
    g_G = \binom{k_{m-1}+q_m}{2}=\binom{k_m}{2},
\end{equation}
and
\begin{equation}\label{eq:canonical}
    (k_m-2)\cdot D_m\sim (k_m-2)(k_m+1)\cdot u_{m,1}\sim  K_{G}.
\end{equation}
We also get that for every integer $-2\leq l\leq k_{m-1}+q_m=k_m$,
\begin{equation}\label{eq:D_mInequality}
    r_G(l\cdot D_m)=r_G\big(l(k_m+1)\cdot u_{m,1}\big)\geq \frac{l(l+3)}{2}.
\end{equation}
The inequality implies that for $1\leq l\leq k_m-2$, we have
$
\gon_{\frac{l(l+3)}{2}}(G)\leq l(k_m-1).
$
Because the gonality sequence is strictly increasing, we can extend this bound to the other gonalities. For every $0\leq h\leq l$,
\begin{equation}\label{eq:upperGonalityBound}
    \gon_{\frac{l(l+3)}{2}-h}(G)\leq \gon_{\frac{l(l+3)}{2}}(G)-h= l(k_m+1)-h
\end{equation}
Next, we derive the matching lower bound. By the inductive hypothesis,
\[
\gon_{\frac{l(l+1)}{2}}(G')=\gon_{\frac{l(l+3)}{2}-l}(G')=l(k_{m-1}+1)-l=lk_{m-1}.
\]
Moreover, for any $1\leq r\leq \frac{(k_{m-1}-2)(k_{m-1}+1)}{2}=g_{G'}-1$, 
\[
r_{G'}\big(\gon_r(G')\cdot v\big)=r_{G'}\big[\big(l(k_{m-1}+1)-h\big)\cdot v\big]\geq r_{G'}\big(l(k_{m-1}+1)\cdot v\big)-h=\frac{l(l+3)}{2}-h=\gon_r(G').
\]
This means that $\gon_r(G')\cdot v$ realises the $r$-th gonality of $G'$. The same conclusion follows for $r\geq g_{G'}$ from \cref{corollary:corollariesOfRiemann-Roch}. Now the assumptions of \cref{lemma:lowerBoundMainLemma} hold for $G'$ and $G$. Applying the lemma with $k=k_{m-1}$ and $n=q_m$, gives that for all $1\leq l\leq k_{m-1}+q_m-2=k_m-2$,
\[
\gon_{\frac{l(l+1)}{2}}(G)\geq l(k_{m-1}+q_m)=lk_m,
\]
or equivalently,
\[
\gon_{\frac{l(l+3)}{2}-l}(G)\geq l(k_m+1)-l.
\]
Again, since the gonality sequence is strictly increasing, for any $0\leq h\leq l$,
\begin{equation}\label{eq:lowerGonalityBound}
    \gon_{\frac{l(l+3)}{2}-h}(G)\geq \gon_{\frac{l(l+3)}{2}-l}(G)+(l-h)\geq l(k_m+1)-h.
\end{equation}

Combining \cref{eq:upperGonalityBound,eq:lowerGonalityBound},
gives the equality:
\[
\gon_{\frac{l(l+3)}{2}-h}(G)=l(k_m+1)-h,\qquad \text{for any $1\leq l\leq k_m-2$, and $0\leq h\leq l$.}
\]
For $r\geq g_G$, \cref{corollary:corollariesOfRiemann-Roch} gives $\gon_r(G)=g_G+r$. Now as
\[
\frac{l(l+3)}{2}-l-1=\frac{(l-1)(l+2)}{2},\qquad\text{and}\qquad \frac{(k_m-2)(k_m+1)}{2}=g_G-1,
\]
the full gonality sequence is
\[
\gon_r(G)=
    \begin{cases}
        l(k_m+1)-h & \text{if } r < g_G, \\
        g+r  & \text{if } r\geq g_G,
    \end{cases}
 \]  
where $1\leq l\leq k-2$ and $0\leq h\leq l$ are uniquely determined integers such that $r=\frac{l(l+3)}{2}-h$, proving Part (2) of the theorem.
In particular,
\[
\gon_1(G)=k_m,\qquad \gon_2(G)=k_m+1,
\]
and $G$ is gonality-tight. Combining this with \cref{eq:genusOfQuasiBananaGraph} proves Part (1).

Next, we show that if $k_m>2$, all divisors realising the second gonality are linearly equivalent. We first verify that $D_m$ realises the second gonality. 
\begin{equation}\label{eq:D_mrealisesGonalityPart1}
    \deg(D_m)=|V(G)|+1=\sum_{i=1}^mq_i+2=k_{m}+1=\gon_2(G)
\end{equation}
and by \cref{eq:D_mInequality} with $l=1$,
\begin{equation}\label{eq:D_mrealisesGonalityPart2}
r_G(D_m)\geq 2,
\end{equation}
as required.
By \cref{claim:chipFiringOnQuasiBananaGraphs},
\[
D_m=u_{0,1}+\sum_{w\in V(G)}w\sim \left(\sum_{i=1}^{m-1} q_i+2\right)\cdot u_{m-1,1} + \sum_{w\in V(K)}w=(k_{m-1}+1)\cdot v + \sum_{i=1}^{q_m} u_{m,i}=:D_m'.
\]
Because there are two or more edges between $v$ and each vertex of $K$, starting a fire on $u_{m-1,1}$ with respect to $D_m'$ burns all of $G$. Therefore, $D_m'$ is $v$-reduced. 

\begin{claim}\label{claim:onlyvReduced}
    The only $v$-reduced divisor on $G$ which realises the second gonality is $D'_m$.
\end{claim}
\begin{proof}
Assume for a contradiction that $D''_m$ is a $v$-reduced divisor on $G$ of degree $k_m+1$ and rank $2$ with $D_m''\neq D_m'$. By \cref{lemma:canRestrictToComponent}
\[
r_{G'}(D_m''|_{G'})\geq 2.
\]
Since $\gon_2(G')=k_{m-1}+1$, we must have that $\deg(D''|_{G'})\geq k_{m-1}+1$, that is
\[
M:=\deg(D''_m|_{K})=\deg (D''_m)-\deg(D''_m|_{G'})\leq k_m+1-(k_{m-1}+1)=q_m.
\]
We consider the following three cases and arrive at a contradiction at each one.

\casetitle[$M = q_m$]{A}

\medskip
\noindent
If \[
    \deg(D''_m|_{K})=M=q_m=V(K), \qquad\text{and}\qquad  D''_m\neq D'_m,
    \] then there must exist a vertex $u\in V(K)$ such that $D''_m(u)=0$.
    Also,
    \begin{equation}\label{eq:chipsOnvA}
        D''_m(v)\leq \deg(D''_m|_{G'}) = \deg(D''_m)-\deg(D''_m|_{K})=k_m+1-q_m=k_{m-1}+1.
    \end{equation}
    Consider the divisor $D''_m-u-v$. It is not effective, as there are $-1$ chips on $u$. Now, start a fire at $u$ with respect to $D''_m-u-v$. From \cref{eq:chipsOnvA}, there are fewer than $k_{m-1}+1$ chips on $v$, so fire spreads to $v$. Since $D''_m$ is $v$-reduced, the entire graph burns down. Therefore, $D''_m-u-v$ is not equivalent to an effective divisor, contradicting the assumption that $D''_m$ has rank $2$.

\casetitle[$1\leq M \leq q_m-1$]{B}

\medskip
\noindent
First note that
    \begin{equation}\label{eq:chipsOnvB}
        D''_m(v)\leq \deg(D''_m|_{G'}) = \deg(D''_m)-\deg(D''_m|_{K})=k_m+1-M=k_{m-1}+1+(q_m-M),
    \end{equation}
    and $V(K)=q_m\geq 2$.
    Order the vertices of $K$ so that
    \[
    D''_m(u_{m,1})\leq D''_m(u_{m,2})\leq\cdots \leq D''_m(u_{m,q_m}).
    \]
    Now as 
    \[
    \sum_{i=1}^{q_m}D(u_{m,i})=\deg(D''_m|_{K})=M\leq q_m-1,
    \]
    we have
    \[
    D''_m(u_{m,1})= D''_m(u_{m,2})=\cdots = D''_m(u_{m,q_m-M})=0,\qquad\text{and}\qquad D''_m(u_{m,q_m-M+1})\leq 1.
    \]
    Let
    \[
    A=\{u_{m,2},u_{m,3},\ldots, u_{m,q_m-M+1} \},\qquad \text{and} \qquad E=u_{m,1}+u_{m,q_m-M+1}.
    \]
    Consider the divisor $D''_m-E$. It has $-1$ chips on $u_{m,1}$ and $0$ or fewer chips on every vertex in $A$. Therefore, starting a fire at $u_{m,1}$ with respect to $D''_m-E$ burns all of $A$.
    The total number of edges from $A$ to $v$ is equal to 
    \[
    |A|(k_{m-1}+1)=(q_m-M+1)(k_{m-1}+1).
    \]
    From \cref{eq:chipsOnvB}, 
    \[
    (D''_m-E)(v) = D''_m(v)\leq k_{m-1}+1+(q_m-M)< (q_m-M+1)(k_{m-1}+1).
    \]
    Hence the fire reaches $v$ and burns all of $G$, contradicting the assumption that $D_m''$ has rank 2.

\casetitle[$M=0$]{C}

\medskip
\noindent
In this case,
    \[
    \deg(D''_m|_{G'}) =\deg(D''_m)=k_m+1= k_{m-1}+q_m+1.
    \]
    Set-firing $V(G')$ produces the divisor
    \[
    F:=D''_m-q_m(k_{m-1}+1)\cdot v + \sum_{w\in V(K)}(k_{m-1}+1)\cdot w.
    \]
    The rank of $F$ is $2$ as $F\sim D''_m$. For any $u\in V(K)$, starting a fire at $v$ with respect to $F-u$ burns $u$ and then all of $K$. Additionally, as $D''_m$ is $v$-reduced and
    \[
    F|_{G'\setminus\{v\}}=D''_m|_{G'\setminus\{v\}},
    \]
    the fire spreads to the rest of $G'$. Therefore $F-u$ is $v$-reduced and has rank at least $1$. Now by \cref{lemma:canRestrictToComponent}, $F|_{G'}$ has rank at least $1$ as a divisor on $G'$. However, since $k_{m-1}+q_m=k_m>2$,
    \[
    \deg(F|_{G'})=\deg(D''_m|_{G'})-q_m(k_{m-1}+1)=(k_{m-1}+q_m+1)-q_m(k_{m-1}+1)<k_{m-1}=\gon_1(G'),
    \]
    a contradiction.
\end{proof}

Two divisors are linearly equivalent if and only if they are equivalent to the same $v$-reduced divisor. As rank is invariant under equivalence, by \cref{claim:onlyvReduced}, every divisor realising the second gonality of $G$ must be equivalent to $D_m'$. This proves Part (3) of the theorem.

For Part (4), let $D$ be any divisor realising the second gonality. 
If $k_m=\gon_1(G)=2$, then the genus of $G$ is $1$, and the result follows by \cref{corollary:corollariesOfRiemann-Roch}. 
Now assume that $k_m>2$.
As all divisors realising the second gonality are equivalent by Part (3), we have $D\sim D'_m\sim D_m$. Therefore, for every integer $1\leq l\leq k_m$,
\[
l\cdot D\sim l\cdot D_m\qquad\text{and}\qquad r_G(l\cdot D)=r_G(l\cdot D_m).
\]
From the computed gonality sequence, notice that, for every $1\leq l\leq k_m$,
\[
\deg(l\cdot D_m)=l\deg(D_m)=l(|V(G)|+1)=l\left(\sum_{i=1}^mq_i+2\right)=l(k_m+1)=\gon_{\frac{l(l+3)}{2}}(G).
\]
Combining this with the inequality in \cref{eq:D_mInequality}, we conclude that $l\cdot D\sim l\cdot D_m$ realises the $\frac{l(l+3)}{2}$-th gonality.
Finally, \cref{eq:canonical} shows that $(k_m-2)\cdot D\sim (k_m-2)\cdot D_m \sim K_G$, completing the proof.

\subsubsection{The algebraic version of \cref{mainTheorem:propertiesOfQBGraphs}}\label{sec:algebraicA}

As mentioned in the introduction,  it was shown in \cite[Lemma 7.2]{FJKO_semigroupOfGonality} that analogous statements to the first two parts of \cref{mainTheorem:propertiesOfQBGraphs} hold in algebraic geometry. We now show that the third part of the conjecture holds as well.

\begin{lemma}\label{lem:algebraicVersionOfConjecture}
    Let $C$ be a smooth curve with first gonality $k$ and second gonality $k+1$ for some integer $k$. Denote by $D$ a divisor that realises the second gonality. For $1\leq l\leq k-2$, a divisor realises the $\frac{l(l+3)}{2}$-the gonality if and only if it is equivalent to $l\cdot D$.
\end{lemma}

\begin{proof}
   The result is an almost immediate consequence of \cite[Theorem 2.1]{Hartshorne_GorensteinAndNoether}, which Hartshorne attributes to Noether.  By \cite[Lemma 7.1]{FJKO_semigroupOfGonality}, we can assume that $C$ is a plane curve of degree $k+1$ and that the divisor obtained as a  line section of $C$ is equivalent to $D$.
   
   Let $Z$ be a divisor realising the $\frac{l(l+3)}{2}$-th gonality. By \cite[Lemma 7.2]{FJKO_semigroupOfGonality}, the degree of $Z$ is $l\cdot (k+1)$. From Part (2.b) of Noether's theorem, the divisor $Z$ is obtained as the scheme-theoretic intersection of $C$ with another curve $C''$ of degree $l$. Therefore $Z$ is equivalent to $l\cdot D$. 
   
   Similarly, since $D$ is a line section of $C$,  the divisor $l\cdot D$ is the intersection of $C$ with a curve of degree $l$, so the ``only if" direction of Noether's theorem implies that the rank of $l\cdot D$ is $\frac{l\cdot(l+3)}{2}$. 
\end{proof}

\subsection{Proof of \cref{mainTheorem:alternativeDefinitonForQBGraphs}}\label{sec:ProofOfThmA}

The following is an immediate  corollary of the proof of \cref{mainTheorem:propertiesOfQBGraphs}.

\begin{corollary}
    Let $G$ be a quasi-banana graph. Then the divisor
    \[
    u_{0,1}+\sum_{w\in V(G)}w
    \]
    realises the second gonality.
\end{corollary}

\begin{proof}
    In the proof of \cref{mainTheorem:propertiesOfQBGraphs}, it is established that for every quasi-banana graph $G$, the divisor
    \[
    D_m=u_{0,1}+\sum_{w\in V(G)}w
    \]
    has degree equal to the second gonality of $G$ (see \cref{eq:D_mrealisesGonalityPart1}) and rank at least $2$ (see \cref{eq:D_mrealisesGonalityPart2}).
\end{proof}

Moreover, \cref{mainTheorem:propertiesOfQBGraphs} implies that every quasi-banana graph is gonality-tight, regardless of the first gonality. This established the forward implication of \cref{mainTheorem:alternativeDefinitonForQBGraphs}. The converse is precisely \cref{proposition:firstDirectionOfMain}. Therefore, the theorem follows.

\subsection{Proof of \cref{mainTheorem:subdividingEdges}}\label{sec:proofOfTheoremC}

Fix an integer $\ell\geq 1$ and let $G$ be a quasi-banana graph. By \cref{mainTheorem:alternativeDefinitonForQBGraphs}, $G$ is gonality-tight, so there exists a divisor $D$ on $G$ with degree $k+1$ and rank $2$. The divisor $D$ as a divisor on $\sigma_\ell (G)$ has rank $2$ by \cite[Corollary 3.4]{HladkyKralNorine_divisorsTropicalCurves}. So $\gon_2 (\sigma_\ell (G))\leq \deg (D) =k+1$.

As the gonality sequence is strictly increasing, it suffices to show that $\gon_1 (\sigma_\ell (G))\geq k$. We require the following lemma.

\begin{lemma}\label{lemma:subdividedEdges}
    Let $H'$ be a graph with $\gon_1 (H')\geq k$ for some integer $k$. Let $H''$ be the graph obtained by attaching a complete graph $K_n$ to $H'$ along $k+1$ edges between $v\in V(H')$ and each vertex of $K_n$. Let $H$ be the graph obtained by subdividing each edge of $H''$ not entirely in $H'$ exactly $\ell\geq 1$ times. Then $\gon_1(H)\geq k+n$.
\end{lemma}
\begin{proof}
    Let $D$ be a divisor on $H$ with degree $k+n-1$. We will show that the rank of $D$ is less than $1$. Without loss of generality, assume that $D$ is $v$-reduced.

    Assume for a contradiction that $r_H(D)<1$.
    By \cref{lemma:canRestrictToComponent}, $D|_{H'}$ has rank at least $1$ as a divisor on $H'$. Since $\gon_1(H')\geq k$, it follows that $\deg(D|_{H'})\geq k$. Therefore, 
    \[
    \deg\left(D|_{H\setminus V(H')}\right)= \deg(D) - \deg (D|_{H'})\leq n-1.
    \]
    There are $n$ vertices in $K_n$, so there must be some vertex $u\in V(K_n)$ such that $D(u)=0$. We aim to show that $D-u$ is $u$-reduced and therefore not equivalent to an effective divisor. It suffices to show that burning $H$ from $u$ with respect to $D-u$ causes $v$ to burn. Then the rest of the whole graph burns down, since we assumed that $D$ is $v$-reduced. 
    
    Suppose that $v$ does not burn down. Let $u_1,\ldots,u_h$ be the vertices of $K_n$ that burn down when starting a fire at $u$. Consider a burnt vertex $u_i$. To prevent the fire from spreading further, there must be at a chip on each unburnt vertex of $K_n$ or on a vertex in a subdivided edge connecting the two of them. Since there are $n-h$ unburnt vertices of $K_n$, this requires $n-h$ chips in total. 
    Similarly, for each of the $k+1$ subdivided edges connecting $u_i$ and $v$, there must be a chip on a vertex in the subdivided edge or on $v$. This requires additional $k+1$ chips.

    As there are $h$ burnt vertices, $\deg(D)\geq h(n-h+k+1)$. Since $1\leq h\leq n$, we have 
    \[
    \deg(D)\geq h(k+1)+h(n-h) = h(k+n+1-h)\geq k+n,
    \]
    contradicting the assumption that $\deg(D)=k+n-1$. Hence, $D-u$ is $u$-reduced and not equivalent to an effective divisor. It follows that $D$ has rank less than $1$. As $D$ was arbitrary, $\gon (H)\geq k+n$, as required.
\end{proof}

We can now prove \cref{mainTheorem:subdividingEdges}.
Write $G=QB_m(q_1,\ldots,q_m)$ for some integers $q_1,\ldots,q_m\geq 1$ and $m\geq0$. We proceed by induction on $m$.

When $m=0$, the only quasi-banana graph to consider is $QB_0$, consisting of a single vertex and no edges. Then $\sigma_\ell (QB_0)=QB_0$, so $\gon_1(G)\geq k$.

Now, assume that $m\geq 1$ and set $G'=QB_{m-1}(q_1,\ldots, q_{m-1})$. By \cref{mainTheorem:alternativeDefinitonForQBGraphs}, the gonality of quasi-banana graphs equals the number of their vertices, so
\[
\gon_1(G') = 1+\sum_{i=1}^{m-1}q_i\quad\text{and}\quad \gon_1(G') = 1+\sum_{i=1}^{m}q_i.
\]
Assume that $\gon_1(\sigma_\ell(G'))\geq \gon_1 (G')$. The graph $\sigma_\ell (G)$ is obtained from $\sigma_\ell(G')$ by attaching the complete graph $K_{q_m}$ and subdividing edges just as in \cref{lemma:subdividedEdges}. Setting $H'=\sigma_\ell (G'), H=\sigma_\ell (G)$ and applying the lemma gives
\[
\gon_1 (\sigma_\ell(G))\geq \gon_1(\sigma_\ell(G'))+q_m\geq \gon_1 (G') +q_m=\gon_1 (G). 
\]
By induction, we conclude that $\gon_1(\sigma_\ell (G))\geq k$ for any quasi-banana graph $G$ with gonality $k$.
As $\gon_2\sigma_\ell (G))\leq k+1$, this proves the first part of the theorem.

Finally, since  the first two gonalities are invariant under equal subdivisions of the edges, it follows from \cite[Theorem 1.3]{DBSW_gonalityCanBeDifferent} that they also coincide with the first two gonalities of the metric graph $\Gamma$. In particular, $\Gamma$ is gonality-tight. 
\qed

\medskip
 
We finish the paper by showing that
subdividing the edges of a quasi-banana graph by unequal amounts does not necessarily preserve gonality-tightness. 
\begin{example}\label{example:unequal subdivision is not gonality tight}
    The graph $G$ depicted below is obtained by subdiving edges of $QB_3(1,1,1)$.
    \begin{figure}[H]
    \centering
    \begin{tikzpicture}[scale=1.3, node/.style={circle, draw, fill=black, inner sep=1.2pt}]

    % Vertices
    \node[node, label=below left:{$v_1$}] (v1) at (0,0) {};
    \node[node] (v2) at (2,0) {};
    \node[below=0.2cm of v2] {$v_2$};
    %\node[node] (v3) at (3.5,0.8) {};
    %\node[node] (v4) at (3.15,-0.8) {};
    \node[node] (v7) at (3,0) {};
    \node[node] (v8) at (4,0) {};
    %\node[node] (v5) at (3.85,-0.8) {};
    \node[node] (v6) at (5,0) {};
    \node[below=0.2cm of v6] {$v_4$};
    \node[node, label=below right:{$v_5$}] (v9) at (7,0) {};
    \node[node] at ($(v2)!0.5!(v6)+(0,0.82455)$) {};
    \node[node, label=below:{$v_3$}] at ($(v2)!0.5!(v6)+(0,-0.82455)$) {};
    
    % Edges
    \draw[bend left=50]  (v1) to (v2);
    \draw[bend right=50] (v1) to (v2);

    %\draw (v2) -- (v3) -- (v6);
    \draw (v2) -- (v7) -- (v8) -- (v6);
    \draw[bend left=65]  (v6) to (v2);
    \draw[bend right=65]  (v6) to (v2);
    %\draw (v2) -- (v4) -- (v5) -- (v6);

    \draw[bend left=65]  (v6) to (v9);
    \draw[bend left=20]   (v6) to (v9);
    \draw[bend right=20]  (v6) to (v9);
    \draw[bend right=65] (v6) to (v9);
    \end{tikzpicture}
    \end{figure}
    
    \noindent We leave it for the reader to to check that the divisors
    \[
    v_1+v_2+v_4+v_5 \quad \text{and} \quad 2v_1+2v_2+v_3+v_5
    \]
    have ranks $1$ and $2$ respectively, and  that divisors of lower degrees have smaller ranks. In particular, the first two gonalities are 4 and 6, so $G$ is not gonality-tight. 
    
    %realise the first and second gonalities respectively. Now, $\gon_1(G)=4$, $\gon_2(G)=6$ and $g(G)=6$. So by \cite[Example 4.3]{ADMYY_GonalitySequencesOfGraphs},
     %the full gonality sequence of $G$ is $4,6,8,9,10,12,13,14,\ldots$. \simun{Once we compute the first two entries, the rest of the sequence can be computed using Riemann-Roch. Do we need to write a proof of this? I don't think there is a clever way to find $\gon_1$ and $\gon_2$ other than computing the ranks of a lot of divisors. }
     %\yoav{I would write down the divisors realising the first two gonalities, and say that we leave the details as an exercise for the reader. I would include the Riemann--Roch argument, especially since it's not immediately clear to me how that follows. I can see that for gonalities 5 and above, but how do 3 and 4 follow?  } \simun{If there is a divisor $D$ with $\deg(D)= 7,r(D)\geq 3$, then by Riemann-Roch, $r(K-D)\geq 1$ and $\deg(K-D)=3$. This contradicts the fact that $\gon_1=4$. Therefore, $\gon_3\geq 8$. As $\gon_5 = 10$, we must have $\gon_3=8$ and $\gon_4=9$. Instead of writing the whole argument out, I referenced \cite{ADMYY_GonalitySequencesOfGraphs}. Let me know whether that's fine.}
\end{example}

% \simun{Should we mention as a corollary that it follows from \cite[Theorem 1.3]{DBSW_gonalityCanBeDifferent} that for a quasi-banana graph $G$, the metric graph $\Gamma(G)$ is gonality-tight as well? We can also say that the complete metric graph $\Gamma (K_n)$ is gonality-tight by \cite{CoolsPanizzut_gonalityComplete}. ($\Gamma(G)$ is the metric graph obtained by assigning length $1$ to each edge in $G$)}

\bibliographystyle{alpha}
\bibliography{Bibliography}

\end{document}